\documentclass[10pt]{article}
\usepackage{amstext}
\usepackage[fleqn]{amsmath}
\usepackage{amsfonts}
\usepackage{amssymb}
\usepackage{amsthm}
\usepackage{newlfont}
\usepackage{graphicx}
\usepackage{sectsty}
\theoremstyle{plain} \theoremstyle{theorem}
\newtheorem{theorem}{Theorem}[section]
\theoremstyle{example}

\theoremstyle{corollary}
\newtheorem{corollary}{Corollary}[section]
\theoremstyle{lemma}
\newtheorem{lemma}{Lemma}[section]
\theoremstyle{proposition}

\theoremstyle{axiom}

\theoremstyle{notation}

\theoremstyle{fact}

\theoremstyle{definition}
\newtheorem{definition}{Definition}[section]
\theoremstyle{remark}
\newtheorem{remark}{Remark}[section]
\numberwithin{equation}{section}

\begin{document}
%---------------------------------------------------------------------------------------------------------------------------------------------------------------------------------------
\title{On incomplete exponential $\;_{r}R_{s}(P,Q,z)$ matrix function}
%---------------------------------------------------------------------------------------------------------------------------------------------------------------------------------------
\author{Ayman Shehata \thanks{%
E-mail: drshehata2006@yahoo.com, aymanshehata@science.aun.edu.eg}\;
, Ghazi S. Khammsh \thanks{%
ghazikhamash@yahoo.com}\;
,Ajay K. Shukla \thanks{%
ajayshukla2@rediffmail.com} and
Shimaa I. Moustafa \thanks{%
E-mail: shimaa1362011@yahoo.com, shimaa$_{-}$m@science.aun.edu.eg} \\
{\small $^{*,\S}$\; Department of Mathematics, Faculty of Science, Assiut University, Assiut 71516, Egypt.}\\
{\small $^{\ddag}$\; Department of Applied Mathematics and Humanities,}\\
{\small Sardar Vallabhbhai National Institute of Technology, Surat-395 007, Gujarat, India.}\\
{\small $^{\dag}$\; Department of Mathematics, Al-Aqsa University, Gaza Strip, Palestine.}}
%---------------------------------------------------------------------------------------------------------------------------------------------------------------------------------------
\date{}
\maketitle{}
%---------------------------------------------------------------------------------------------------------------------------------------------------------------------------------------
%-------------------------------------------------------------------------------------------------------------------------------------------------------------------------------------
\begin{abstract}
%-------------------------------------------------------------------------------------------------------------------------------------------------------------------------------------
The recurrence matrix relations, differentiation formulas, and analytical and fractional integral properties of incomplete gamma matrix functions  $\gamma(Q, x)$ and $\Gamma(Q, x)$ are all covered in this article.  The generalized incomplete exponential matrix functions with their integral representations functions have been examined, along with some relevant characteristics of these functions such as integral representations  functions . Additionally, the infinite summation relations and formulas for two sequences are shown, along with the generalized incomplete exponential matrix functions with the integral representation, addition formula for addition of two arguments, multiplication formula for multiplication of two arguments, and recurrence matrix relation.
%---------------------------------------------------------------------------------------------------------------------------------------------------------------------------------------
\end{abstract}
%-------------------------------------------------------------------------------------------------------------------------------------------------------------------------------------
%-------------------------------------------------------------------------------------------------------------------------------------------------------------------------------------

%---------------------------------------------------------------------------------------------------------------------------------------------------------------------------------------
\section{Introduction}
%-------------------------------------------------------------------------------------------------------------------------------------------------------------------------------------
In 1950, the Tricomi studied the theory of incomplete gamma functions \cite{tr}. These functions are very essential special functions
and that are used in number of problems in mathematical, physical, astrophysics, applied statistics, probability theory, statistical sciences, approximation theory, and engineering. These
functions are also useful in the study of various transform such as Laplace transforms, Fourier transform and probability
theory \cite{ba, emo}.

In \cite{cz, et, cq, ssj}, the researchers investigated the generalized incomplete Gamma functions $\Gamma(\nu;x)$ and $\gamma(\nu;x)$ and their applications. By using each of these classical incomplete Gamma functions$\Gamma(\nu;x)$ and $\gamma(\nu;x)$, the incomplete Pochhammer symbols $(\nu;x)_{n}$ and $[\nu;x]_{n}$ were defined, for $x\geq0$ and $\nu,n\in \mathcal{C}$, by Srivastava et al. \cite{sca}. Srivastava et al. \cite{sca} established many properties of the classical incomplete hypergeometric functions by the help of these representations. Later, the incomplete special functions were established by researchers work \cite{sca} which have a particularly strong connection with fractional calculus are the incomplete gamma and incomplete beta functions in \cite{sca}.

Throughout this paper, for a matrix $A$ in $\Bbb{C}^{N\times N}$, its spectrum $\sigma(A)$ denotes the set of all eigenvalues of $A$. The two-norm will be denoted by $||A||_{2}$ and it is defined by (see \cite{j3, j4})
%-----------------------------------------------------------------------------------------------------------------------------------------------------------------
\begin{equation*}
\begin{split}
||A||_{2}=\sup_{x\neq 0}\frac{||Ax||_{2}}{||x||_{2}},
\end{split}
\end{equation*}
%-----------------------------------------------------------------------------------------------------------------------------------------------------------------
where for a vector $x$ in $\Bbb{C}^{N}$, $||x||_2=(x^Tx)^\frac{1}{2}$ is the Euclidean norm of $x$.
%-----------------------------------------------------------------------------------------------------------------------------------------------------------------

Let us denote the real numbers $M(A)$ and $m(A)$ as in the following
%-----------------------------------------------------------------------------------------------------------------------------------------------------------------
\begin{equation}
\begin{split}
\mathbb{M}(A)=\max\{Re(z): z\in \sigma(A)\};\quad \mathbf{m}(A)=\min\{Re(z): z\in \sigma(A)\}.\label{1.1}
\end{split}
\end{equation}
%-----------------------------------------------------------------------------------------------------------------------------------------------------------------
If $\Phi(z)$ and $\Psi(z)$ are holomorphic functions of the complex variable $z$, which are defined in an open set $\Omega$ of the complex plane, and $A$, $B$
are matrices in $\Bbb{C}^{N\times N}$ with $\sigma(A)\subset\Omega$ and $\sigma(B) \subset \Omega$, such that $AB=BA$, then from the properties of the matrix functional calculus in \cite{nd}, it follows that
%-----------------------------------------------------------------------------------------------------------------------------------------------------------------
\begin{equation*}
\begin{split}
\Phi(A)\Psi(B)=\Psi(B)\Phi(A).
\end{split}
\end{equation*}
%-------------------------------------------------------------------------------------------------------------------------------------------------------------------------------------
Throughout this study, a matrix polynomial of degree $n$ in $x$ means an expression of the form
%-------------------------------------------------------------------------------------------------------------------------------------------------------------------------------------
\begin{equation*}
\begin{split}
P_{n}(x)=A_{n}x^{n}+A_{n-1}x^{n-1}+\ldots+A_{1}x+A_{0},
\end{split}
\end{equation*}
%-------------------------------------------------------------------------------------------------------------------------------------------------------------------------------------
where $x$ is a real variable or complex variable, $A_j$, for $0<j<n$ and $A_{n}\neq \textbf{0}$, are complex matrices in $\Bbb{C}^{N\times N}$.
%-----------------------------------------------------------------------------------------------------------------------------------------------------------------

We recall that the reciprocal gamma function denoted by $\Gamma^{-1}(z)=\frac{1}{\Gamma(z)}$ is an entire function of the complex variable $z$ and thus for any matrix $A$ in $\Bbb{C}^{N\times N}$, $\Gamma^{-1}(A)$ is a well defined matrix. Furthermore, if $A$ is a matrix such that
%-----------------------------------------------------------------------------------------------------------------------------------------------------------------
\begin{equation}
\begin{split}
A+nI\quad \text{is\; an invertible matrix\; for\; all\; integers}\; n\geq0,\label{1.2}
\end{split}
\end{equation}
%-----------------------------------------------------------------------------------------------------------------------------------------------------------------
where $I$ is the identity matrix in $\Bbb{C}^{N\times N}$, then from \cite{j2}, it follows that
%-----------------------------------------------------------------------------------------------------------------------------------------------------------------
\begin{equation}
\begin{split}
(A)_{n}=A(A+I)\ldots(A+(n-1)I)=\Gamma{(A+nI)}\Gamma^{-1}{(A)}\;;\quad n\geq 1\;;\quad (A)_{0}=I.\label{1.3}
\end{split}
\end{equation}
%-----------------------------------------------------------------------------------------------------------------------------------------------------------------
If $k$ is large enough so then for $k>\|B\|$, then we will mention to the following relation
which existed in J\'{o}dar and Cort\'{e}s \cite{j3, j4} in the form
%-----------------------------------------------------------------------------------------------------------------------------------------------------------------
\begin{equation*}
\begin{split}
\|(B+kI)^{-1}\|\leq \frac{1}{k-\|B\|}\;;\quad k>\|B\|.
\end{split}
\end{equation*}
%-----------------------------------------------------------------------------------------------------------------------------------------------------------------
If $\Theta(k,n)$ and $\Omega(k,n)$ are matrices in $\Bbb{C}^{N\times N}$ for $n\geq 0$, $k\geq 0$, then in an analogous way to the proof of Lemma \textbf{11} \cite{j3}, it follows that
%-----------------------------------------------------------------------------------------------------------------------------------------------------------------
\begin{equation}
\begin{split}
\sum_{n=0}^{\infty}\sum_{\ell=0}^{\infty}\Theta(k,n)=\sum_{n=0}^{\infty}\sum_{\ell=0}^{[\frac{1}{2}n]}\Theta(k,n-2k),\\
\sum_{n=0}^{\infty}\sum_{\ell=0}^{\infty}\Omega(k,n)=\sum_{n=0}^{\infty}\sum_{\ell=0}^{n}\Omega(k,n-k).\label{1.4}
\end{split}
\end{equation}
%-----------------------------------------------------------------------------------------------------------------------------------------------------------------
Similarly to (\ref{1.4}), we can write
%-----------------------------------------------------------------------------------------------------------------------------------------------------------------
\begin{equation}
\begin{split}
\sum_{n=0}^{\infty}\sum_{\ell=0}^{[\frac{1}{2}n]}\Theta(k,n)=\sum_{n=0}^{\infty}\sum_{\ell=0}^{\infty}\Theta(k,n+2k),\\
\sum_{n=0}^{\infty}\sum_{\ell=0}^{n}\Omega(k,n)=\sum_{n=0}^{\infty}\sum_{\ell=0}^{\infty}\Omega(k,n+k).\label{1.5}
\end{split}
\end{equation}
%-----------------------------------------------------------------------------------------------------------------------------------------------------------------
The hypergeometric matrix function $_{2}F_{1}(A,B;C;z)$ has been given in the form (see \cite{j3, j4})
%-----------------------------------------------------------------------------------------------------------------------------------------------------------------
\begin{equation}
\begin{split}
_{2}F_{1}(A,B;C;z)=\sum_{\ell=0}^{\infty}\frac{(A)_k(B)_k[(C)_k]^{-1}}{n!}z^k,\label{1.6}
\end{split}
\end{equation}
%-----------------------------------------------------------------------------------------------------------------------------------------------------------------
for matrices $A$, $B$ and $C$ in $\Bbb{C}^{N\times N}$ such that $C+nI$ is an invertible matrix for all integers $n\geq 0$ and for $|z|<1$. It has been
seen by J\'{o}dar and Cort\'{e}s \cite{j3} that the series is absolutely convergent for $|z|=1$ when
%-----------------------------------------------------------------------------------------------------------------------------------------------------------------
\begin{equation*}
\begin{split}
\mathbf{m}(C)> \mathbb{M}(A)+\mathbb{M}(B),
\end{split}
\end{equation*}
%-----------------------------------------------------------------------------------------------------------------------------------------------------------------
where $\mathbf{m}(A)$ and $\mathbb{M}(A)$ in (\ref{1.1}) for any matrix $A$ in $\Bbb{C}^{N\times N}$.
%-------------------------------------------------------------------------------------------------------------------------------------------------------------------------------------
\begin{definition}
%-------------------------------------------------------------------------------------------------------------------------------------------------------------------------------------
For $p$ and $q$ are finite positive integers, the generalized hypergeometric matrix function defined as (see \cite{sa1})
%-----------------------------------------------------------------------------------------------------------------------------------------------------------------
\begin{equation}
\begin{split}
&\;_{p}F_{q}(A_{1},A_{2},\ldots,A_{p};B_{1},B_{2},\ldots,B_{q};z)\\
=&\sum_{k=0}^{\infty}\frac{z^{k}}{k!}(A_{1})_{k}(A_{2})_{k}\ldots(A_{p})_{k}[(B_{1})_{k}]^{-1}[(B_{2})_{k}]^{-1}\ldots[(B_{q})_{k}]^{-1}\\
=&\sum_{k=0}^{\infty}\frac{z^{k}}{k!}\prod_{i=1}^{p}(A_{i})_{k}\bigg{[}\prod_{j=1}^{q}(B_{i})_{k}\bigg{]}^{-1},\label{1.7}
\end{split}
\end{equation}
%-----------------------------------------------------------------------------------------------------------------------------------------------------------------
where $A_{i}$ ; $1\leq i\leq p$ and $B_{j}$; $1\leq j\leq q$ are matrices in $\Bbb{C}^{N\times N}$ such that
%-----------------------------------------------------------------------------------------------------------------------------------------------------------------
\begin{equation}
\begin{split}
B_{j}+nI \quad \text{are\; invertible\; matrices\; for\; all\; integers}\; n\geq0,\label{1.8}
\end{split}
\end{equation}
%-----------------------------------------------------------------------------------------------------------------------------------------------------------------
%-----------------------------------------------------------------------------------------------------------------------------------------------------------------
\begin{enumerate}
  \item If $p\leq q$, then the power series (\ref{1.7}) converges for all finite $z$.
  \item If $p=q+1$, then the power series (\ref{1.7}) is convergent for $|z|<1$ and diverges for $|z|>1$.
  \item If $p>q+1$, then the power series (\ref{1.7}) diverges for all $z$, $z\neq0$.
\end{enumerate}
%-----------------------------------------------------------------------------------------------------------------------------------------------------------------
%-----------------------------------------------------------------------------------------------------------------------------------------------------------------
%-----------------------------------------------------------------------------------------------------------------------------------------------------------------
\begin{enumerate}
%-----------------------------------------------------------------------------------------------------------------------------------------------------------------
\item If $p=q+1$, then the power series (\ref{1.7}) is absolutely convergent for $|z|=1$ when
%-----------------------------------------------------------------------------------------------------------------------------------------------------------------
\begin{equation}
\begin{split}
\sum_{j=1}^{q}\mathbf{m}(B_{j})>\sum_{i=1}^{p}\mathbb{M}(A_{j}).\label{1.9}
\end{split}
\end{equation}
%-------------------------------------------------------------------------------------------------------------------------------------------------------------------------------------
\item If $p=q+1$, then the power series (\ref{1.7}) is conditionally convergent for $|z|=1$ when
%-----------------------------------------------------------------------------------------------------------------------------------------------------------------
\begin{equation}
\begin{split}
\sum_{i=0}^{p}\mathbb{M}(A_{j})-1<\sum_{j=0}^{q}\mathbf{m}(B_{j})\leq\sum_{i=0}^{p}\mathbb{M}(A_{j}).\label{1.10}
\end{split}
\end{equation}
%-----------------------------------------------------------------------------------------------------------------------------------------------------------------
\item If $p=q+1$, then the power series (\ref{1.7}) is diverges for $|z|=1$ when
%-----------------------------------------------------------------------------------------------------------------------------------------------------------------
\begin{equation}
\begin{split}
\sum_{j=0}^{q}\mathbf{m}(B_{j})\leq\sum_{i=0}^{p}\mathbb{M}(A_{j})-1,\label{1.11}
\end{split}
\end{equation}
%-----------------------------------------------------------------------------------------------------------------------------------------------------------------
where $\mathbb{M}(A_{i})$ and $\mathbf{m}(B_{j})$ are defined in (\ref{1.1}).
%-----------------------------------------------------------------------------------------------------------------------------------------------------------------
\end{enumerate}
%-------------------------------------------------------------------------------------------------------------------------------------------------------------------------------------
%-------------------------------------------------------------------------------------------------------------------------------------------------------------------------------------
\end{definition}
%-------------------------------------------------------------------------------------------------------------------------------------------------------------------------------------
\begin{definition}
%-----------------------------------------------------------------------------------------------------------------------------------------------------------------
Let us suppose that $P$, $Q$, $Re(P)>0$, $Re(Q)>0$, $A_{i}$; $Re(A_{i})>0$, $1\leq i\leq r$ and $B_{j}$; $Re(B_{j})>0$, $1\leq j\leq s$ are matrices in $\Bbb{C}^{N\times N}$ such that
%-----------------------------------------------------------------------------------------------------------------------------------------------------------------
\begin{equation}
\begin{split}
B_{j}+\ell I \quad \text{are\; invertible\; matrices\; for\; all\; integers}\; \ell \geq0,\label{1.12}
\end{split}
\end{equation}
%-----------------------------------------------------------------------------------------------------------------------------------------------------------------
where $r$ and $s$ are finite positive integers. Then we define the $\;_{r}R_{s}(P,Q,z)$ matrix function as
%-----------------------------------------------------------------------------------------------------------------------------------------------------------------
\begin{equation}
\begin{split}
&\;_{r}R_{s}(A_{1},A_{2},\ldots,A_{r};B_{1},B_{2},\ldots,B_{s};P,Q;z)\\
=&\sum_{\ell=0}^{\infty}\frac{z^{\ell}}{\ell !}(A_{1})_{\ell}(A_{2})_{\ell}\ldots(A_{r})_{\ell}[(B_{1})_{\ell}]^{-1}[(B_{2})_{\ell}]^{-1}\ldots[(B_{s})_{\ell}]^{-1}\Gamma^{-1}(\ell P+Q)\\
=&\sum_{\ell=0}^{\infty}\frac{z^{\ell}}{\ell !}\prod_{i=1}^{r}(A_{i})_{\ell}\bigg{[}\prod_{j=1}^{s}(B_{j})_{\ell}\bigg{]}^{-1}\Gamma^{-1}(\ell P+Q)=\sum_{\ell=0}^{\infty}W_{\ell},\label{1.13}
\end{split}
\end{equation}
%---------------------------------------------------------------------------------------------------------------------------------------------------------------------------------------
\end{definition}
%-----------------------------------------------------------------------------------------------------------------------------------------------------------------
where $W_{\ell}=\frac{z^{\ell}}{\ell !}\prod_{i=1}^{r}(A_{i})_{\ell}\bigg{[}\prod_{j=1}^{s}(B_{j})_{\ell}\bigg{]}^{-1}\Gamma^{-1}(\ell P+Q)$.
%-----------------------------------------------------------------------------------------------------------------------------------------------------------------

%-----------------------------------------------------------------------------------------------------------------------------------------------------------------
The last limit shows that
%-----------------------------------------------------------------------------------------------------------------------------------------------------------------
\begin{enumerate}
  \item If $r\leq s+1$, then the power series in (\ref{1.13}) converges for all finite $z$.
  \item If $r=s+2$, then the power series in (\ref{1.13}) converges for all $|z|<1$ and diverges for all $|z|>1$.
  \item If $r>s+2$, then the power series in (\ref{1.13}) diverges for $z\neq0$.
\end{enumerate}
%-----------------------------------------------------------------------------------------------------------------------------------------------------------------
In analogous to theorem 3 in \cite{j3}, we can state the following:
%-----------------------------------------------------------------------------------------------------------------------------------------------------------------
\begin{theorem}
%-----------------------------------------------------------------------------------------------------------------------------------------------------------------
\begin{enumerate}
%-----------------------------------------------------------------------------------------------------------------------------------------------------------------
\item If $r=s+2$, then the power series in (\ref{1.13}) is absolutely convergent on the circle $|z|=1$ when
%-----------------------------------------------------------------------------------------------------------------------------------------------------------------
\begin{equation}
\begin{split}
\sum_{j=1}^{s}m(B_{j})-\sum_{i=1}^{r}M(A_{i})>0.\label{1.14}
\end{split}
\end{equation}
%-----------------------------------------------------------------------------------------------------------------------------------------------------------------
\item If $r=s+2$, then the power series (\ref{1.13}) is conditionally convergent for $|z|=1$ when
%-----------------------------------------------------------------------------------------------------------------------------------------------------------------
\begin{equation}
\begin{split}
\sum_{i=0}^{r}M(A_{i})-1<\sum_{j=0}^{s}m(B_{j})\leq\sum_{i=0}^{p}M(A_{i}).\label{1.15}
\end{split}
\end{equation}
%-----------------------------------------------------------------------------------------------------------------------------------------------------------------
\item If $r=s+2$, then the power series (\ref{1.13}) is diverges for $|z|=1$ when
%-----------------------------------------------------------------------------------------------------------------------------------------------------------------
\begin{equation}
\begin{split}
\sum_{j=0}^{s}m(B_{j})\leq\sum_{i=0}^{r}M(A_{i})-1\label{1.16}
\end{split}
\end{equation}
%-----------------------------------------------------------------------------------------------------------------------------------------------------------------
where $M(A_{i})$; $1\leq i\leq r$  and $m(B_{j})$; $1\leq j\leq s$ are defined in (\ref{1.1}).
%---------------------------------------------------------------------------------------------------------------------------------------------------------------------------------------------
\end{enumerate}
%---------------------------------------------------------------------------------------------------------------------------------------------------------------------------------------------
\end{theorem}
%---------------------------------------------------------------------------------------------------------------------------------------------------------------------------------------------
Thus $\;_{r}R_{s}$ is an entire function of $z$ when $\|P+I\|>0$.
%---------------------------------------------------------------------------------------------------------------------------------------------------------------------------------------------
\begin{remark} Let $A_{i}$\,; $1\leq i\leq r$ and $B_{j}$\,; $1\leq j\leq s$ be matrices in $\Bbb{C}^{N\times N}$ satisfying (\ref{1.12}) and all matrices are commutative. For $P=Q=A_{1}=I$ in (\ref{1.13}) reduces to
%---------------------------------------------------------------------------------------------------------------------------------------------------------------------------------------------
\begin{equation}
\begin{split}
&\;_{r}R_{s}(I,A_{2},\ldots,A_{p};B_{1},B_{2},\ldots,B_{s};I,I;z)\\
=&\sum_{\ell=0}^{\infty}\frac{z^{\ell}}{k!}\prod_{i=2}^{p}(A_{i})_{\ell}\bigg{[}\prod_{j=1}^{s}(B_{j})_{\ell}\bigg{]}^{-1}\Gamma^{-1}(kP+Q)=\sum_{\ell=0}^{\infty}W_{\ell}\\
=&\;_{r-1}F_{s}(A_{2},\ldots,A_{p};B_{1},B_{2},\ldots,B_{s};z),\label{1.17}
\end{split}
\end{equation}
%---------------------------------------------------------------------------------------------------------------------------------------------------------------------------------------------
where $\;_{r-1}F_{s}$ is generalized hypergeometric matrix function in (\ref{1.7}).
%---------------------------------------------------------------------------------------------------------------------------------------------------------------------------------------------
%-----------------------------------------------------------------------------------------------------------------------------------------------------------------
%---------------------------------------------------------------------------------------------------------------------------------------------------------------------------------------------
\end{remark}
%---------------------------------------------------------------------------------------------------------------------------------------------------------------------------------------------
Numerous researcher such as Shukla and Prajapati [17], Salim [14]and Sarma [15] studied incomplete exponential functions and gave its generalization. The motive of the current work is to investigation the analytical and fractional integral properties of new type incomplete exponential matrix functions. This function is an amalgamation of generalized incomplete exponential function which plays an important role in the theory of mathematical analysis, fractional calculus and statistics and has significant applications in the field of free electron laser equations and fractional kinetic equations. For literature survey of fractional integral operators, researchers can refer to the papers of Riesz (1949), Srivastava and Tomovski (2009) and Gupta et al. (2019) and books of Miller and Ross (1993), Milovanovic´ and Rassias (2014), Govil et al. (2017), Gupta et al. (2018) and Gupta and Rassias (2019).

The organization of the our present paper is as follows: In Section 2, we discuss the incomplete gamma matrix functions $\gamma(Q, x)$ and $\Gamma(Q, x)$ and establish several properties on the partial derivatives of these matrix functions. In Section 3, we define a new or known extension of the type incomplete exponential matrix functions and investigate certain properties of these matrix functions. Finally, some concluding remarks and consequences of the results are exhibited discussed in Section 4.
%-------------------------------------------------------------------------------------------------------------------------------------------------------------------------------------
%-------------------------------------------------------------------------------------------------------------------------------------------------------------------------------------
%---------------------------------------------------------------------------------------------------------------------------------------------------------------------------------------------
\section{Incomplete gamma matrix functions $\gamma(Q, x)$ and $\Gamma(Q, x)$}
%---------------------------------------------------------------------------------------------------------------------------------------------------------------------------------------------
In this section, we discuss the our main properties of incomplete gamma matrix functions $\gamma(Q, x)$ and $\Gamma(Q, x)$.
%-------------------------------------------------------------------------------------------------------------------------------------------------------------------------------------
\begin{definition}
%-------------------------------------------------------------------------------------------------------------------------------------------------------------------------------------
%---------------------------------------------------------------------------------------------------------------------------------------------------------------------------------------------
Let $Q$ be a matrix in$\Bbb{C}^{N\times N}$ such that $\check{\mu}(Q)>0$, and $x$ be a positive real number. Then the familiar incomplete gamma matrix functions $\gamma(Q, x)$ and $\Gamma(Q, x)$ are defined as \cite{sj}
%---------------------------------------------------------------------------------------------------------------------------------------------------------------------------------------
\begin{equation}
\begin{split}
\gamma(Q,x) := \int_{0}^{x} u^{Q-I} e^{-u} du, \;\;\;\;\;\;\; (\Re(Q)>0; x\geq0),\label{2.1}
\end{split}
\end{equation}
%---------------------------------------------------------------------------------------------------------------------------------------------------------------------------------------------
and
%---------------------------------------------------------------------------------------------------------------------------------------------------------------------------------------------
\begin{equation}
\begin{split}
\Gamma(Q,x) := \int_{x}^{\infty} u^{Q-I} e^{-u} du, \;\;\;\;\;\;\; (\Re(Q)>0; x\geq0).\label{2.2}
\end{split}
\end{equation}
%-------------------------------------------------------------------------------------------------------------------------------------------------------------------------------------
\end{definition}
%-------------------------------------------------------------------------------------------------------------------------------------------------------------------------------------
%---------------------------------------------------------------------------------------------------------------------------------------------------------------------------------------------
The incomplete gamma matrix functions $\gamma(Q, x)$ and $\Gamma(Q, x)$ satisfy the following relation:
%---------------------------------------------------------------------------------------------------------------------------------------------------------------------------------------------
\begin{equation}
\begin{split}
\gamma(Q,x) + \Gamma(Q,x) = \Gamma(Q), \;\;\;\; (\Re(Q)>0).\label{2.3}
\end{split}
\end{equation}
%---------------------------------------------------------------------------------------------------------------------------------------------------------------------------------------
\begin{theorem}
%---------------------------------------------------------------------------------------------------------------------------------------------------------------------------------------
The $\gamma(Q, x)$ and $\Gamma(Q, x)$ satisfies the recurrence matrix relations
%---------------------------------------------------------------------------------------------------------------------------------------------------------------------------------------
\begin{equation}
\begin{split}
\gamma(Q+I,x)=Q\gamma(Q,x)-e^{-x}x^{Q}\label{2.4}
\end{split}
\end{equation}
%---------------------------------------------------------------------------------------------------------------------------------------------------------------------------------------
and
%---------------------------------------------------------------------------------------------------------------------------------------------------------------------------------------
\begin{equation}
\begin{split}
\Gamma(Q+I,x)=Q\Gamma(Q,x)+e^{-x}x^{Q}.\label{2.5}
\end{split}
\end{equation}
%---------------------------------------------------------------------------------------------------------------------------------------------------------------------------------------
%---------------------------------------------------------------------------------------------------------------------------------------------------------------------------------------
\end{theorem}
%---------------------------------------------------------------------------------------------------------------------------------------------------------------------------------------
%---------------------------------------------------------------------------------------------------------------------------------------------------------------------------------------
\begin{proof}
%---------------------------------------------------------------------------------------------------------------------------------------------------------------------------------------
Using (\ref{2.1}), we have %---------------------------------------------------------------------------------------------------------------------------------------------------------------------------------------
\begin{equation*}
\begin{split}
\gamma(Q+I,x)= \int_{0}^{x} t^{Q+I-I} e^{-t} dt=\int_{0}^{x} t^{Q} e^{-t} dt,
\end{split}
\end{equation*}
%---------------------------------------------------------------------------------------------------------------------------------------------------------------------------------------------
$u= t^{Q}\Longrightarrow du=Qt^{Q-I}dt$ and $dv=e^{-t}dt\Longrightarrow v=-e^{-t}$
%---------------------------------------------------------------------------------------------------------------------------------------------------------------------------------------
\begin{equation*}
\begin{split}
\gamma(Q+I,x)=-t^{Q}e^{-t}|_{0}^{x}+Q\int_{0}^{x} t^{Q-I} e^{-t} dt\\
=-x^{Q}e^{-x}+Q\gamma(Q,x)
\end{split}
\end{equation*}
%---------------------------------------------------------------------------------------------------------------------------------------------------------------------------------------------
%---------------------------------------------------------------------------------------------------------------------------------------------------------------------------------------------
\begin{equation*}
\begin{split}
\Gamma(Q+I,x) = \int_{x}^{\infty} t^{Q+I-I} e^{-t} dt.
\end{split}
\end{equation*}
%---------------------------------------------------------------------------------------------------------------------------------------------------------------------------------------------
%---------------------------------------------------------------------------------------------------------------------------------------------------------------------------------------------
$u= t^{Q}\Longrightarrow du=Qt^{Q-I}dt$ and $dv=e^{-t}dt\Longrightarrow v=-e^{-t}$
%---------------------------------------------------------------------------------------------------------------------------------------------------------------------------------------
%---------------------------------------------------------------------------------------------------------------------------------------------------------------------------------------------
\begin{equation*}
\begin{split}
\Gamma(Q+I,x) = \int_{x}^{\infty} t^{Q+I-I} e^{-t} dt= \int_{x}^{\infty} t^{Q} e^{-t} dt\\
=-t^{Q}e^{-t}|_{x}^{0}+Q\int_{x}^{\infty} t^{Q-I} e^{-t} dt\\
=x^{Q}e^{-x}+Q\Gamma(Q,x).
\end{split}
\end{equation*}
%---------------------------------------------------------------------------------------------------------------------------------------------------------------------------------------------
%---------------------------------------------------------------------------------------------------------------------------------------------------------------------------------------
\end{proof}
%---------------------------------------------------------------------------------------------------------------------------------------------------------------------------------------
%---------------------------------------------------------------------------------------------------------------------------------------------------------------------------------------
\begin{theorem}
%---------------------------------------------------------------------------------------------------------------------------------------------------------------------------------------
The matrix functions $\gamma(Q, x)$ and $\Gamma(Q, x)$ satisfy the recurrence relation
%---------------------------------------------------------------------------------------------------------------------------------------------------------------------------------------
\begin{equation}
\begin{split}
\gamma(Q+2I,x)-(Q+(x+1)I)\gamma(Q+I,x)+Qx\gamma(Q,x)=\mathbf{0}\label{2.6}
\end{split}
\end{equation}
%---------------------------------------------------------------------------------------------------------------------------------------------------------------------------------------
and
%---------------------------------------------------------------------------------------------------------------------------------------------------------------------------------------
\begin{equation}
\begin{split}
\Gamma(Q+2I,x)-(Q+(x+1)I)\Gamma(Q+I,x)+Qx\Gamma(Q,x)=\mathbf{0}.\label{2.7}
\end{split}
\end{equation}
%---------------------------------------------------------------------------------------------------------------------------------------------------------------------------------------
%---------------------------------------------------------------------------------------------------------------------------------------------------------------------------------------
\end{theorem}
%---------------------------------------------------------------------------------------------------------------------------------------------------------------------------------------
%---------------------------------------------------------------------------------------------------------------------------------------------------------------------------------------
%-------------------------------------------------------------------------------------------------------------------------------------------------------------------------------------
\begin{definition}
%-------------------------------------------------------------------------------------------------------------------------------------------------------------------------------------
%---------------------------------------------------------------------------------------------------------------------------------------------------------------------------------------------
%---------------------------------------------------------------------------------------------------------------------------------------------------------------------------------------------
The incomplete Pochhammer symbols $(Q;x)_n$ and $[Q;x]_n$ can also be represented in terms of incomplete gamma matrix functions $\gamma(Q,x)$ and $\Gamma(Q,x)$ as, (Normalized matrix functions are defined:)
%---------------------------------------------------------------------------------------------------------------------------------------------------------------------------------------------
\begin{equation}
\begin{split}
(Q;x)_n := \gamma(Q+nI,x)\Gamma^{-1}(Q), \;\;\;\;\; (Q, n \in \mathbf{C}; x\geq0)\label{2.8}
\end{split}
\end{equation}
%---------------------------------------------------------------------------------------------------------------------------------------------------------------------------------------------
and
%---------------------------------------------------------------------------------------------------------------------------------------------------------------------------------------------
\begin{equation}
\begin{split}
[Q;x]_n := \Gamma(Q+nI,x)\Gamma^{-1}(Q), \;\;\;\;\; (Q, n \in \mathbf{C}; x\geq0).\label{2.9}
\end{split}
\end{equation}
%-------------------------------------------------------------------------------------------------------------------------------------------------------------------------------------
\end{definition}
%-------------------------------------------------------------------------------------------------------------------------------------------------------------------------------------
%---------------------------------------------------------------------------------------------------------------------------------------------------------------------------------------------
%---------------------------------------------------------------------------------------------------------------------------------------------------------------------------------------------
In the view of (\ref{2.8}) and (\ref{2.9}) give the following decomposition relation:
%---------------------------------------------------------------------------------------------------------------------------------------------------------------------------------------------
\begin{equation}
\begin{split}
(Q;x)_n + [Q;x]_n = (Q)_n &:=  \Gamma(Q+n)\Gamma^{-1}(Q) = \begin{cases}
     I \;\;\;\;\;\;\;\;\;\;\;\;\;\; \hspace{2cm}          (n=0; Q\in\mathbf{C}\setminus\{0\})\\
      Q(Q+I)\dots(Q+(n-1)I) \;\;  \hspace{0.4cm}(n\in\mathbf{N}; Q\in\mathbf{C}).\label{2.10}
   \end{cases}
\end{split}
\end{equation}
%---------------------------------------------------------------------------------------------------------------------------------------------------------------------------------------------
where $(Q)_n$ is the Pochhammer symbol (\ref{1.3}). In particular, the generalized pochhammer symbol $(Q)_{kn}$ can be represented in the following form \cite{sa1}:
%---------------------------------------------------------------------------------------------------------------------------------------------------------------------------------------------
\begin{equation}
\begin{split}
(Q)_{kn} = k^{kn} \ {\left( {\frac{Q }{k}} \right)_n}{\left( {\frac{{Q  + I}}{k}} \right)_n} \cdots  \cdots {\left( {\frac{{Q  + (k - 1)I}}{k}} \right)_n}.\label{2.11}
\end{split}
\end{equation}
%---------------------------------------------------------------------------------------------------------------------------------------------------------------------------------------------
The matrix functions (\ref{2.1}) and (\ref{2.2}) are closely connected with the particular case $P=I$ of the confluent hypergeometric matrix functions $\Phi(P;Q;x)$ and $U(P;Q;x)$ via
%---------------------------------------------------------------------------------------------------------------------------------------------------------------------------------------
\begin{equation}
\begin{split}
\gamma(Q,x) = Q^{-1}x^{Q}e^{-x}\Phi(I;Q+I;x),\\
\gamma(Q,x) = Q^{-1}\Phi(Q;Q+I;x)\label{2.12}
\end{split}
\end{equation}
%---------------------------------------------------------------------------------------------------------------------------------------------------------------------------------------------
and
%---------------------------------------------------------------------------------------------------------------------------------------------------------------------------------------
\begin{equation}
\begin{split}
\Gamma(Q,x) = x^{Q}e^{-x}U(I;Q+I;x),\\
\Gamma(Q,x) = e^{-x}U(I-Q;I-Q;x).\label{2.13}
\end{split}
\end{equation}
%---------------------------------------------------------------------------------------------------------------------------------------------------------------------------------------------
%---------------------------------------------------------------------------------------------------------------------------------------------------------------------------------------
%---------------------------------------------------------------------------------------------------------------------------------------------------------------------------------------
\begin{theorem}
%---------------------------------------------------------------------------------------------------------------------------------------------------------------------------------------
%---------------------------------------------------------------------------------------------------------------------------------------------------------------------------------------
The matrix functions $\gamma(Q, x)$ and $\Gamma(Q, x)$ satisfy the differentiation formulas
%---------------------------------------------------------------------------------------------------------------------------------------------------------------------------------------
%---------------------------------------------------------------------------------------------------------------------------------------------------------------------------------------
\begin{equation}
\begin{split}
\frac{\partial\gamma(Q,x)}{\partial x} =x^{Q-I}e^{-x}\label{2.14}
\end{split}
\end{equation}
%---------------------------------------------------------------------------------------------------------------------------------------------------------------------------------------
and
%---------------------------------------------------------------------------------------------------------------------------------------------------------------------------------------
\begin{equation}
\begin{split}
\frac{\partial\Gamma(Q,x)}{\partial x} =-x^{Q-I}e^{-x}.\label{2.15}
\end{split}
\end{equation}
%---------------------------------------------------------------------------------------------------------------------------------------------------------------------------------------
\end{theorem}
%---------------------------------------------------------------------------------------------------------------------------------------------------------------------------------------
%---------------------------------------------------------------------------------------------------------------------------------------------------------------------------------------
\begin{proof}
%---------------------------------------------------------------------------------------------------------------------------------------------------------------------------------------
Differentiate the following integral with respect to the parameter $t$ in the form
%---------------------------------------------------------------------------------------------------------------------------------------------------------------------------------------
\begin{equation}
\begin{split}
u(x,t)=\int_{\phi(t)}^{\psi(t)}f(x,t)dx,\\
\frac{du(x,t)}{dt}=\int_{\phi(t)}^{\psi(t)}\frac{\partial f(x,t)}{\partial t}dx+f(\psi(t),t)\frac{d\psi(t)}{dt}-f(\phi(t),t)\frac{d\phi(t)}{dt}.\label{2.16}
\end{split}
\end{equation}
%---------------------------------------------------------------------------------------------------------------------------------------------------------------------------------------
%---------------------------------------------------------------------------------------------------------------------------------------------------------------------------------------
%---------------------------------------------------------------------------------------------------------------------------------------------------------------------------------------
Using the equation above, the derivative of the upper incomplete gamma matrix function with respect to the parameter $x$ are
%---------------------------------------------------------------------------------------------------------------------------------------------------------------------------------------
\begin{equation*}
\begin{split}
\frac{\partial\gamma(Q,x)}{\partial x} =x^{Q-I}e^{-x}
\end{split}
\end{equation*}
%---------------------------------------------------------------------------------------------------------------------------------------------------------------------------------------
and
%---------------------------------------------------------------------------------------------------------------------------------------------------------------------------------------
\begin{equation*}
\begin{split}
\frac{\partial\Gamma(Q,x)}{\partial x} =-x^{Q-I}e^{-x}.
\end{split}
\end{equation*}
%---------------------------------------------------------------------------------------------------------------------------------------------------------------------------------------
%---------------------------------------------------------------------------------------------------------------------------------------------------------------------------------------
\end{proof}
%---------------------------------------------------------------------------------------------------------------------------------------------------------------------------------------
%---------------------------------------------------------------------------------------------------------------------------------------------------------------------------------------
\begin{theorem}
%---------------------------------------------------------------------------------------------------------------------------------------------------------------------------------------
The matrix functions $\gamma(Q, x)$ and $\Gamma(Q, x)$ satisfy the differential equations
%---------------------------------------------------------------------------------------------------------------------------------------------------------------------------------------
\begin{equation}
\begin{split}
\frac{d^{2}\gamma(Q,x)}{d x^{2}}+\bigg{[}I+\frac{I-Q}{x}\bigg{]}\frac{d \gamma(Q,x)}{d x}=\mathbf{0}\label{2.17}
\end{split}
\end{equation}
%---------------------------------------------------------------------------------------------------------------------------------------------------------------------------------------
and
%---------------------------------------------------------------------------------------------------------------------------------------------------------------------------------------
\begin{equation}
\begin{split}
\frac{d^{2}\Gamma(Q,x)}{d x^{2}}+\bigg{[}I+\frac{I-Q}{x}\bigg{]}\frac{d \Gamma(Q,x)}{d x}=\mathbf{0}.\label{2.18}
\end{split}
\end{equation}
%---------------------------------------------------------------------------------------------------------------------------------------------------------------------------------------
%---------------------------------------------------------------------------------------------------------------------------------------------------------------------------------------
%---------------------------------------------------------------------------------------------------------------------------------------------------------------------------------------
\end{theorem}
%---------------------------------------------------------------------------------------------------------------------------------------------------------------------------------------
%---------------------------------------------------------------------------------------------------------------------------------------------------------------------------------------
\begin{proof}
%---------------------------------------------------------------------------------------------------------------------------------------------------------------------------------------
Successively (\ref{2.17}) and (\ref{2.18}) can be proved in an analogous manner as (a similar way as) in Theorem 2.3.
%---------------------------------------------------------------------------------------------------------------------------------------------------------------------------------------
%---------------------------------------------------------------------------------------------------------------------------------------------------------------------------------------
\end{proof}
%---------------------------------------------------------------------------------------------------------------------------------------------------------------------------------------
\begin{theorem}
%---------------------------------------------------------------------------------------------------------------------------------------------------------------------------------------
%---------------------------------------------------------------------------------------------------------------------------------------------------------------------------------------
The matrix functions $\gamma(Q, x)$ and $\Gamma(Q, x)$ satisfies the following identities
%---------------------------------------------------------------------------------------------------------------------------------------------------------------------------------------
%---------------------------------------------------------------------------------------------------------------------------------------------------------------------------------------
\begin{equation}
\begin{split}
\frac{d^{n}}{dx^{n}}\bigg{[}x^{-Q}\Gamma(Q,x)\bigg{]} =(-1)^{n}x^{-Q-nI}\Gamma(Q+nI,x),\label{2.19}
\end{split}
\end{equation}
%---------------------------------------------------------------------------------------------------------------------------------------------------------------------------------------
%---------------------------------------------------------------------------------------------------------------------------------------------------------------------------------------
\begin{equation}
\begin{split}
\frac{d^{n}}{dx^{n}}\bigg{[}x^{-Q}\gamma(Q,x)\bigg{]} =(-1)^{n}x^{-Q-nI}\gamma(Q+nI,x),\label{2.20}
\end{split}
\end{equation}
%---------------------------------------------------------------------------------------------------------------------------------------------------------------------------------------
%---------------------------------------------------------------------------------------------------------------------------------------------------------------------------------------
\begin{equation}
\begin{split}
\frac{d^{n}}{dx^{n}}\bigg{[}e^{x}\Gamma(Q,x)\bigg{]} =(-1)^{n}(I-Q)_{n}e^{x}\Gamma(Q-nI,x)\label{2.21}
\end{split}
\end{equation}
%---------------------------------------------------------------------------------------------------------------------------------------------------------------------------------------
and
%---------------------------------------------------------------------------------------------------------------------------------------------------------------------------------------
\begin{equation}
\begin{split}
\frac{d^{n}}{dx^{n}}\bigg{[}e^{x}\gamma(Q,x)\bigg{]} =(-1)^{n}(I-Q)_{n}e^{x}\gamma(Q-nI,x).\label{2.22}
\end{split}
\end{equation}
%---------------------------------------------------------------------------------------------------------------------------------------------------------------------------------------
%---------------------------------------------------------------------------------------------------------------------------------------------------------------------------------------
%---------------------------------------------------------------------------------------------------------------------------------------------------------------------------------------
\end{theorem}
%---------------------------------------------------------------------------------------------------------------------------------------------------------------------------------------

%---------------------------------------------------------------------------------------------------------------------------------------------------------------------------------------
%---------------------------------------------------------------------------------------------------------------------------------------------------------------------------------------
\begin{proof}
%---------------------------------------------------------------------------------------------------------------------------------------------------------------------------------------
Respectively (\ref{2.19})-(\ref{2.22}) can be proved in a similar way as in Theorem 3.2 2.3.%---------------------------------------------------------------------------------------------------------------------------------------------------------------------------------------
\end{proof}
%---------------------------------------------------------------------------------------------------------------------------------------------------------------------------------------

%---------------------------------------------------------------------------------------------------------------------------------------------------------------------------------------
\begin{theorem}
%---------------------------------------------------------------------------------------------------------------------------------------------------------------------------------------
The matrix function $\Gamma(Q, x)$ satisfy the integral representation
%---------------------------------------------------------------------------------------------------------------------------------------------------------------------------------------
%---------------------------------------------------------------------------------------------------------------------------------------------------------------------------------------
\begin{equation}
\begin{split}
\Gamma(Q,x) =2\int_{\sqrt{x}}^{\infty}u^{2Q-I}e^{-u^{2}}du\label{2.23}
\end{split}
\end{equation}
%---------------------------------------------------------------------------------------------------------------------------------------------------------------------------------------
and
%---------------------------------------------------------------------------------------------------------------------------------------------------------------------------------------
\begin{equation}
\begin{split}
\gamma(Q,x) =e^{-x}x^{Q}\int_{0}^{1}(1-u)^{Q-I}e^{ux}du.\label{2.24}
\end{split}
\end{equation}
%---------------------------------------------------------------------------------------------------------------------------------------------------------------------------------------
\end{theorem}
%---------------------------------------------------------------------------------------------------------------------------------------------------------------------------------------
\begin{proof}
%---------------------------------------------------------------------------------------------------------------------------------------------------------------------------------------
The substitution $t=u^{2}$ in (\ref{2.2}), we obtain (\ref{2.23}). The transformation $t=x(1-u)$ in (\ref{2.1}), we give (\ref{2.23}).
%---------------------------------------------------------------------------------------------------------------------------------------------------------------------------------------
\end{proof}
%---------------------------------------------------------------------------------------------------------------------------------------------------------------------------------------
%---------------------------------------------------------------------------------------------------------------------------------------------------------------------------------------
\begin{theorem}
%---------------------------------------------------------------------------------------------------------------------------------------------------------------------------------------
The following indefinite integrals are readily obtained using integration by parts
%---------------------------------------------------------------------------------------------------------------------------------------------------------------------------------------
\begin{equation}
\begin{split}
\int x^{A-I}\Gamma(Q,x)dx=A^{-1}\bigg{[}x^{A}\Gamma(Q,x)-\Gamma(Q+A,x)\bigg{]}+c \label{2.25}
\end{split}
\end{equation}
%---------------------------------------------------------------------------------------------------------------------------------------------------------------------------------------
and
%---------------------------------------------------------------------------------------------------------------------------------------------------------------------------------------
\begin{equation}
\begin{split}
\int x^{A-I}\gamma(Q,x)dx=A^{-1}\bigg{[}x^{A}\gamma(Q,x)-\gamma(Q+A,x)\bigg{]}+c.\label{2.26}
\end{split}
\end{equation}
%---------------------------------------------------------------------------------------------------------------------------------------------------------------------------------------
%---------------------------------------------------------------------------------------------------------------------------------------------------------------------------------------
\end{theorem}
%---------------------------------------------------------------------------------------------------------------------------------------------------------------------------------------
%---------------------------------------------------------------------------------------------------------------------------------------------------------------------------------------
\begin{theorem}
%---------------------------------------------------------------------------------------------------------------------------------------------------------------------------------------
Assume that $Q$ is a matrix in $\Bbb{C}^{N\times N}$ such that $\check{\mu}(Q)>0$, and $x$ be a positive real number. Then, we have the
following Integral Representation
%---------------------------------------------------------------------------------------------------------------------------------------------------------------------------------------
%---------------------------------------------------------------------------------------------------------------------------------------------------------------------------------------
\begin{equation}
\begin{split}
\int_{0}^{t} \Gamma\bigg{(}Q,\frac{1}{x}\bigg{)}dx=t\Gamma\bigg{(}Q,\frac{1}{t}\bigg{)}-\Gamma\bigg{(}Q-I,\frac{1}{t}\bigg{)}.\label{2.27}
\end{split}
\end{equation}
%---------------------------------------------------------------------------------------------------------------------------------------------------------------------------------------
%---------------------------------------------------------------------------------------------------------------------------------------------------------------------------------------
\end{theorem}
%---------------------------------------------------------------------------------------------------------------------------------------------------------------------------------------
%---------------------------------------------------------------------------------------------------------------------------------------------------------------------------------------
\begin{theorem}
%---------------------------------------------------------------------------------------------------------------------------------------------------------------------------------------
Assume that $Q$ is a matrix in $\Bbb{C}^{N\times N}$ such that $\check{\mu}(Q)>0$, and $x$ be a positive real number. Then, we have the
following Integral formulas
%---------------------------------------------------------------------------------------------------------------------------------------------------------------------------------------
\begin{equation}
\begin{split}
\int_{0}^{\infty}e^{-Px}\gamma(Q,x)dx =P^{-1}\gamma(Q)(I+P)^{-Q};Re(P)>0,Re(Q)>-1,\label{2.28}
\end{split}
\end{equation}
%---------------------------------------------------------------------------------------------------------------------------------------------------------------------------------------
%---------------------------------------------------------------------------------------------------------------------------------------------------------------------------------------
\begin{equation}
\begin{split}
\int_{0}^{\infty}e^{-Px}\Gamma(Q,x)dx =P^{-1}\Gamma(Q)[I-(I+P)^{-Q}];Re(P)>-1,Re(Q)>-1,\label{2.29}
\end{split}
\end{equation}
%---------------------------------------------------------------------------------------------------------------------------------------------------------------------------------------
\begin{equation}
\begin{split}
\int_{0}^{\infty}x^{P-I}\gamma(Q,x)dx =-P^{-1}\Gamma(Q+P);Re(P)<0,Re(Q+P)>0\label{2.30}
\end{split}
\end{equation}
%---------------------------------------------------------------------------------------------------------------------------------------------------------------------------------------
and
%---------------------------------------------------------------------------------------------------------------------------------------------------------------------------------------
%---------------------------------------------------------------------------------------------------------------------------------------------------------------------------------------
\begin{equation}
\begin{split}
\int_{0}^{\infty}x^{P-I}\Gamma(Q,x)dx =P^{-1}\Gamma(Q+P);Re(P)>0,Re(Q+P)>0.\label{2.31}
\end{split}
\end{equation}
%---------------------------------------------------------------------------------------------------------------------------------------------------------------------------------------
%---------------------------------------------------------------------------------------------------------------------------------------------------------------------------------------
\end{theorem}
%---------------------------------------------------------------------------------------------------------------------------------------------------------------------------------------

%---------------------------------------------------------------------------------------------------------------------------------------------------------------------------------------
\begin{theorem}
%---------------------------------------------------------------------------------------------------------------------------------------------------------------------------------------
The integral representation for the matrix function $\gamma*(Q,x)$ satisfies the series
%---------------------------------------------------------------------------------------------------------------------------------------------------------------------------------------
\begin{equation}
\begin{split}
\gamma(Q+nI,x) =(Q)_{n}\gamma(Q,x)-x^{Q}e^{-x}\sum_{r=0}^{n-1}\Gamma(Q+nI)\Gamma^{-1}(Q+(r+1)I)x^{r}\label{2.32}
\end{split}
\end{equation}
%---------------------------------------------------------------------------------------------------------------------------------------------------------------------------------------
and
%---------------------------------------------------------------------------------------------------------------------------------------------------------------------------------------
%---------------------------------------------------------------------------------------------------------------------------------------------------------------------------------------
\begin{equation}
\begin{split}
\Gamma(Q+nI,x) =(Q)_{n}\Gamma(Q,x)+x^{Q}e^{-x}\sum_{r=0}^{n-1}\Gamma(Q+nI)\Gamma^{-1}(Q+(r+1)I)x^{r}.\label{2.33}
\end{split}
\end{equation}
%---------------------------------------------------------------------------------------------------------------------------------------------------------------------------------------
%---------------------------------------------------------------------------------------------------------------------------------------------------------------------------------------
\end{theorem}
%---------------------------------------------------------------------------------------------------------------------------------------------------------------------------------------
%---------------------------------------------------------------------------------------------------------------------------------------------------------------------------------------
%-------------------------------------------------------------------------------------------------------------------------------------------------------------------------------------
\begin{definition}
%-------------------------------------------------------------------------------------------------------------------------------------------------------------------------------------
%---------------------------------------------------------------------------------------------------------------------------------------------------------------------------------------------
Now, we define the incomplete gamma matrix functions $\gamma*(Q,x)$ in the form
%---------------------------------------------------------------------------------------------------------------------------------------------------------------------------------------
\begin{equation}
\begin{split}
\gamma*(Q,x) = x^{-Q}\gamma(Q,x)\Gamma^{-1}(Q),\label{2.34}
\end{split}
\end{equation}
%---------------------------------------------------------------------------------------------------------------------------------------------------------------------------------------------
which is an entire matrix function in a as well as in $\Bbb{C}^{N\times N}$ and real valued for real $Re(Q)>0$ and
real $x$.
%---------------------------------------------------------------------------------------------------------------------------------------------------------------------------------------
%-------------------------------------------------------------------------------------------------------------------------------------------------------------------------------------
\end{definition}
%-------------------------------------------------------------------------------------------------------------------------------------------------------------------------------------
In general, the matrix function (\ref{2.1}) is expressed in terms of the confluent hypergeometric matrix function to give ([1], [124])
%---------------------------------------------------------------------------------------------------------------------------------------------------------------------------------------
\begin{equation}
\begin{split}
\gamma*(Q,x) = \Gamma^{-1}(Q+I)\Phi(Q;Q+I;x).\label{2.35}
\end{split}
\end{equation}
%---------------------------------------------------------------------------------------------------------------------------------------------------------------------------------------------
%---------------------------------------------------------------------------------------------------------------------------------------------------------------------------------------
\begin{theorem}
%---------------------------------------------------------------------------------------------------------------------------------------------------------------------------------------
The integral representation for the matrix function $\gamma*(Q,x)$ is
%---------------------------------------------------------------------------------------------------------------------------------------------------------------------------------------
\begin{equation}
\begin{split}
\gamma*(Q,x)=\Gamma^{-1}(Q)\int_{0}^{1} t^{Q-I}e^{-xt}dt.\label{2.36}
\end{split}
\end{equation}
%---------------------------------------------------------------------------------------------------------------------------------------------------------------------------------------
%---------------------------------------------------------------------------------------------------------------------------------------------------------------------------------------
\end{theorem}
%---------------------------------------------------------------------------------------------------------------------------------------------------------------------------------------

%---------------------------------------------------------------------------------------------------------------------------------------------------------------------------------------
\begin{theorem}
%---------------------------------------------------------------------------------------------------------------------------------------------------------------------------------------
The recurrence formula for the matrix function $\gamma*(Q,x)$ is
%---------------------------------------------------------------------------------------------------------------------------------------------------------------------------------------
\begin{equation}
\begin{split}
x\gamma*(Q+I,x) =\gamma*(Q,x)-e^{-x}\Gamma^{-1}(Q+I).\label{2.37}
\end{split}
\end{equation}
%---------------------------------------------------------------------------------------------------------------------------------------------------------------------------------------
\end{theorem}
%---------------------------------------------------------------------------------------------------------------------------------------------------------------------------------------
%---------------------------------------------------------------------------------------------------------------------------------------------------------------------------------------
\begin{lemma}
%---------------------------------------------------------------------------------------------------------------------------------------------------------------------------------------------
The following limit for function $\gamma*(Q,x)$ as
%---------------------------------------------------------------------------------------------------------------------------------------------------------------------------------------
\begin{equation}
\begin{split}
\lim_{x\rightarrow \infty}x^{Q}\gamma*(Q,x) =I.\label{2.38}
\end{split}
\end{equation}
%---------------------------------------------------------------------------------------------------------------------------------------------------------------------------------------
%---------------------------------------------------------------------------------------------------------------------------------------------------------------------------------------
\end{lemma}
%---------------------------------------------------------------------------------------------------------------------------------------------------------------------------------------
\begin{theorem}
%---------------------------------------------------------------------------------------------------------------------------------------------------------------------------------------
The series representations of the matrix function $\gamma*(Q,x)$ is
%---------------------------------------------------------------------------------------------------------------------------------------------------------------------------------------
\begin{equation}
\begin{split}
\gamma*(Q,x) =e^{-x}\sum_{k=0}^{\infty}x^{k}\Gamma^{-1}(Q+(k+1)I)\label{2.39}
\end{split}
\end{equation}
%---------------------------------------------------------------------------------------------------------------------------------------------------------------------------------------
and
%---------------------------------------------------------------------------------------------------------------------------------------------------------------------------------------
\begin{equation}
\begin{split}
\gamma*(Q,x)=\Gamma^{-1}(Q)\sum_{k=0}^{\infty}\frac{(-1)^{k}x^{k}}{k!}(Q+kI)^{-1}.\label{2.40}
\end{split}
\end{equation}
%---------------------------------------------------------------------------------------------------------------------------------------------------------------------------------------
%---------------------------------------------------------------------------------------------------------------------------------------------------------------------------------------
\end{theorem}
%---------------------------------------------------------------------------------------------------------------------------------------------------------------------------------------
\begin{proof}
%---------------------------------------------------------------------------------------------------------------------------------------------------------------------------------------
%---------------------------------------------------------------------------------------------------------------------------------------------------------------------------------------
From (\ref{2.24}), we get
%---------------------------------------------------------------------------------------------------------------------------------------------------------------------------------------
\begin{equation}
\begin{split}
\gamma(Q,x) =e^{-x}x^{Q}\sum_{k=0}^{\infty}x^{k}[(Q)_{k+1}]^{-1}\label{2.41}
\end{split}
\end{equation}
%---------------------------------------------------------------------------------------------------------------------------------------------------------------------------------------
where $Re(Q)\neq0,-1,-2, \ldots$.
%---------------------------------------------------------------------------------------------------------------------------------------------------------------------------------------
From (\ref{2.34}) and (\ref{2.41}), we obtain (\ref{2.39}).
%---------------------------------------------------------------------------------------------------------------------------------------------------------------------------------------

The substitution of the series representation of $e^{-u}$ in (\ref{2.1}) yields
%---------------------------------------------------------------------------------------------------------------------------------------------------------------------------------------
\begin{equation}
\begin{split}
\gamma(Q,x) =x^{Q}\sum_{k=0}^{\infty}\frac{(-1)^{k}x^{k}}{k!}(Q+kI)^{-1}.\label{2.42}
\end{split}
\end{equation}
%---------------------------------------------------------------------------------------------------------------------------------------------------------------------------------------
%---------------------------------------------------------------------------------------------------------------------------------------------------------------------------------------
From (\ref{2.32}) and (\ref{2.42}), we obtain (\ref{2.40}).
%---------------------------------------------------------------------------------------------------------------------------------------------------------------------------------------
%---------------------------------------------------------------------------------------------------------------------------------------------------------------------------------------
\end{proof}
%---------------------------------------------------------------------------------------------------------------------------------------------------------------------------------------
\begin{theorem}
%---------------------------------------------------------------------------------------------------------------------------------------------------------------------------------------
The matrix function $\gamma^{*}(Q,x)$ satisfy the following property
%---------------------------------------------------------------------------------------------------------------------------------------------------------------------------------------
%---------------------------------------------------------------------------------------------------------------------------------------------------------------------------------------
\begin{equation}
\begin{split}
\frac{d^{n}}{dx^{n}}\bigg{[}x^{Q}e^{x}\gamma^{*}(Q,x)\bigg{]} =x^{Q-nI}e^{x}\gamma^{*}(Q-nI,x).\label{2.43}
\end{split}
\end{equation}
%---------------------------------------------------------------------------------------------------------------------------------------------------------------------------------------
%---------------------------------------------------------------------------------------------------------------------------------------------------------------------------------------
%---------------------------------------------------------------------------------------------------------------------------------------------------------------------------------------
\end{theorem}
%---------------------------------------------------------------------------------------------------------------------------------------------------------------------------------------
%---------------------------------------------------------------------------------------------------------------------------------------------------------------------------------------
%---------------------------------------------------------------------------------------------------------------------------------------------------------------------------------------
\begin{theorem}
%---------------------------------------------------------------------------------------------------------------------------------------------------------------------------------------
The matrix functions $\Gamma(Q, x)$ and $\gamma^{*}(Q,x)$ satisfy the following properties
%---------------------------------------------------------------------------------------------------------------------------------------------------------------------------------------
\begin{equation}
\begin{split}
\frac{d^{2}Y}{dx^{2}}-\bigg{[}I+\frac{I-Q}{x}\bigg{]}\frac{dY}{dx}+\frac{I-Q}{x^{2}}Y=0, Y=e^{x}x^{I-Q}\Gamma(Q,x)\label{2.44}
\end{split}
\end{equation}
%---------------------------------------------------------------------------------------------------------------------------------------------------------------------------------------
and
%---------------------------------------------------------------------------------------------------------------------------------------------------------------------------------------
\begin{equation}
\begin{split}
x\frac{d^{2}\gamma^{*}(Q,x)}{dx^{2}}+\bigg{[}Q+(x+1)I\bigg{]}\frac{d\gamma^{*}(Q,x)}{dx}+Q\gamma^{*}(Q,x)=0.\label{2.45}
\end{split}
\end{equation}
%---------------------------------------------------------------------------------------------------------------------------------------------------------------------------------------
\end{theorem}
%---------------------------------------------------------------------------------------------------------------------------------------------------------------------------------------------
\section{On incomplete exponential matrix function}
%---------------------------------------------------------------------------------------------------------------------------------------------------------------------------------------
In this section, we introduce and study several further properties of the incomplete exponential matrix functions.
%---------------------------------------------------------------------------------------------------------------------------------------------------------------------------------------
\begin{definition}
%---------------------------------------------------------------------------------------------------------------------------------------------------------------------------------------
Let $Q$ be a matrix in $\Bbb{C}^{N\times N}$ such that $\check{\mu}(Q)>0$, and $x$ be a positive real number. We define the incomplete exponential matrix functions as follows,
%---------------------------------------------------------------------------------------------------------------------------------------------------------------------------------------
%---------------------------------------------------------------------------------------------------------------------------------------------------------------------------------------------
\begin{equation}
\begin{split}
e\left( {(x,u);Q } \right): = \sum\limits_{n = 0}^\infty  \gamma \left( Q  + nI,x \right)\Gamma^{-1}\left( Q  + nI \right)\frac{{{u^n}}}{{n!}}\label{3.1}
\end{split}
\end{equation}
%---------------------------------------------------------------------------------------------------------------------------------------------------------------------------------------------
and
%---------------------------------------------------------------------------------------------------------------------------------------------------------------------------------------------
\begin{equation}
\begin{split}
E\left( {(x,u);Q } \right): = \sum\limits_{n = 0}^\infty  \Gamma \left( Q  + nI,x \right)\Gamma^{-1}\left( Q  + nI \right)\frac{{{u^n}}}{{n!}}.\label{3.2}
\end{split}
\end{equation}
%---------------------------------------------------------------------------------------------------------------------------------------------------------------------------------------
\end{definition}
%---------------------------------------------------------------------------------------------------------------------------------------------------------------------------------------

%---------------------------------------------------------------------------------------------------------------------------------------------------------------------------------------------
From (\ref{3.1}) and (\ref{3.2}), we obtain the following result,
%---------------------------------------------------------------------------------------------------------------------------------------------------------------------------------------------
\begin{equation}
\begin{split}
e\left( {(x,u);Q } \right) + E\left( {(x,u);Q } \right) = {e^{uI}}.\label{3.3}
\end{split}
\end{equation}
%---------------------------------------------------------------------------------------------------------------------------------------------------------------------------------------------
%---------------------------------------------------------------------------------------------------------------------------------------------------------------------------------------
\begin{lemma}
%---------------------------------------------------------------------------------------------------------------------------------------------------------------------------------------
The following integral representations holds true for the incomplete exponential matrix functions $e\left( {(x,u);Q } \right)$ and $E\left( {(x,u);Q } \right)$:
%---------------------------------------------------------------------------------------------------------------------------------------------------------------------------------------
\begin{equation}
\begin{split}
e\left( {(x,u);Q } \right) &= \Gamma^{-1}(Q) \int_0^x t^{Q-1} \ e^{-t} \left(\sum\limits_{n = 0}^\infty [(Q)_n]^{-1} \frac{(ut)^n}{n!} \right) dt \\ &=\Gamma^{-1}(Q) \int_0^x t^{Q-1} e^{-t}  {}_0F_1(-;Q;tu) dt\label{3.4}
\end{split}
\end{equation}
%---------------------------------------------------------------------------------------------------------------------------------------------------------------------------------------
and
%---------------------------------------------------------------------------------------------------------------------------------------------------------------------------------------------
\begin{equation}
\begin{split}
E\left( {(x,u);Q } \right) &= \Gamma^{-1}(Q) \int_x^\infty t^{Q-1} \ e^{-t} \left(\sum\limits_{n = 0}^\infty [(Q)_n]^{-1} \frac{(ut)^n}{n!} \right) dt \\ &= \Gamma^{-1}(Q) \int_x^\infty t^{Q-1}  e^{-t} {}_0F_1(-;Q;tu) dt.\label{3.5}
\end{split}
\end{equation}
%---------------------------------------------------------------------------------------------------------------------------------------------------------------------------------------
%---------------------------------------------------------------------------------------------------------------------------------------------------------------------------------------
\end{lemma}
%---------------------------------------------------------------------------------------------------------------------------------------------------------------------------------------
In sequel to the study of the generalized incomplete exponential matrix functions (\ref{3.1}) and (\ref{3.2}), we define the following generalized incomplete exponential matrix functions as,
%---------------------------------------------------------------------------------------------------------------------------------------------------------------------------------------
%---------------------------------------------------------------------------------------------------------------------------------------------------------------------------------------------
\begin{equation}
\begin{split}
&\;_{r}e_{s}=\;_{r}e_{s}(x,P,Q;z) = {\displaystyle \;_{r}e_{s}\!\left[(x,P,Q;z)\left |{\begin{matrix} A_1,A_2,\dots,A_r\\ B_1,B_2,\dots,B_s\end{matrix}} \right.\right]},\\
=&\sum_{\ell=0}^{\infty}\frac{z^{\ell}}{\ell !}\prod_{i=1}^{r}(A_{i})_{\ell}\bigg{[}\prod_{j=1}^{s}(B_{j})_{\ell}\bigg{]}^{-1}\gamma(\ell P+Q,x)\Gamma^{-1}(\ell P+Q)\\
=&\sum_{\ell=0}^{\infty}\Psi_{\ell}(z),\label{3.6}
\end{split}
\end{equation}
%---------------------------------------------------------------------------------------------------------------------------------------------------------------------------------------
and
%---------------------------------------------------------------------------------------------------------------------------------------------------------------------------------------------
\begin{equation}
\begin{split}
&\;_{r}E_{s}=\;_{r}E_{s}(x,P,Q;z) = {\displaystyle \;_{r}E_{s}\!\left[(x,P,Q;z)\left |{\begin{matrix} A_1,A_2,\dots,A_r\\ B_1,B_2,\dots,B_s\end{matrix}} \right.\right]},\\
=&\sum_{\ell=0}^{\infty}\frac{z^{\ell}}{\ell !}\prod_{i=1}^{r}(A_{i})_{\ell}\bigg{[}\prod_{j=1}^{s}(B_{j})_{\ell}\bigg{]}^{-1}\Gamma(\ell P+Q,x)\Gamma^{-1}(\ell P+Q)\\
=&\sum_{\ell=0}^{\infty}\Omega_{\ell}(z),\label{3.7}
\end{split}
\end{equation}
%---------------------------------------------------------------------------------------------------------------------------------------------------------------------------------------

where $P,Q \in \mathbf{C}$, \ $A_1, A_2,\dots,A_r \in \mathbf{C}$ and $B_1, B_2,\dots,B_s \in \mathbf{C}\setminus\mathbf{Z}_0^-$ provided that the series defined by (\ref{3.6}) and (\ref{3.7}) on r.h.s. converges.

From (\ref{3.6}) and (\ref{3.7}), we obtain the following decomposition formula:
%---------------------------------------------------------------------------------------------------------------------------------------------------------------------------------------
%---------------------------------------------------------------------------------------------------------------------------------------------------------------------------------------------
\begin{equation}
\begin{split}
{\displaystyle \;_{r}e_{s}\!\left[(x,P,Q;z)\left |{\begin{matrix} A_1,A_2,\dots,A_r\\ B_1,B_2,\dots,B_s\end{matrix}} \right.\right]} + {\displaystyle \;_{r}E_{s}\!\left[(x,P,Q;z)\left |{\begin{matrix} A_1,A_2,\dots,A_r\\ B_1,B_2,\dots,B_s\end{matrix}} \right.\right]} = {}_{r}F_{s}\left[{\left.{\begin{matrix} A_1,A_2,\dots,A_r\\ B_1,B_2,\dots,B_s\end{matrix}} \right|z}\right],\label{3.8}
\end{split}
\end{equation}
%---------------------------------------------------------------------------------------------------------------------------------------------------------------------------------------
%---------------------------------------------------------------------------------------------------------------------------------------------------------------------------------------
where ${}_{r}F_{s}(\cdot)$ recognize as the generalized hypergeometric matrix function defined by (\ref{1.7}).
%---------------------------------------------------------------------------------------------------------------------------------------------------------------------------------------
\begin{remark}
%---------------------------------------------------------------------------------------------------------------------------------------------------------------------------------------
For $r=s=0$ and $P=I$, (\ref{3.6}) and (\ref{3.7}) reduces to incomplete exponential matrix functions (\ref{3.1}) and (\ref{3.2}):
%---------------------------------------------------------------------------------------------------------------------------------------------------------------------------------------
%---------------------------------------------------------------------------------------------------------------------------------------------------------------------------------------------
\begin{equation}
\begin{split}
_0e_0(x,I,Q;z) &= {\displaystyle _0e_0\!\left[(x,I,Q;z)\left |{\begin{matrix} -\\ -\end{matrix}} \right.\right]},\\
&= \sum_{\ell=0}^{\infty} \gamma( Q+\ell I,x)\Gamma^{-1}( Q+\ell I) \frac{z^\ell}{\ell !} = e\left( {(x,z);Q} \right),\label{3.9}
\end{split}
\end{equation}
%---------------------------------------------------------------------------------------------------------------------------------------------------------------------------------------
and
%---------------------------------------------------------------------------------------------------------------------------------------------------------------------------------------------
\begin{equation}
\begin{split}
_0E_0(x,1,Q;z) &= {\displaystyle _0E_0\!\left[(x,I,Q;z)\left |{\begin{matrix} -\\ -\end{matrix}} \right.\right]},\\
&= \sum_{\ell=0}^{\infty} \Gamma( Q+\ell I,x)\Gamma^{-1}( Q+\ell I) \frac{z^\ell}{\ell !} = E\left( {(x,z);Q } \right).\label{3.10}
\end{split}
\end{equation}
%---------------------------------------------------------------------------------------------------------------------------------------------------------------------------------------

%---------------------------------------------------------------------------------------------------------------------------------------------------------------------------------------
\end{remark}
%---------------------------------------------------------------------------------------------------------------------------------------------------------------------------------------
\begin{theorem}
%---------------------------------------------------------------------------------------------------------------------------------------------------------------------------------------
%-----------------------------------------------------------------------------------------------------------------------------------------------------------------
For all integers $n\ge 1$, we have
%-----------------------------------------------------------------------------------------------------------------------------------------------------------------
\begin{equation}
\begin{split}
\;_{r}e_{s}(A_{i}+nI)&=\prod_{k=1}^{n}\bigg{(}A_{i}+(k-1)I\bigg{)}^{-1}\sum_{\ell=0}^{\infty}\prod_{k=1}^{n}(A_{i}+(\ell+k-1)I)\Psi_{\ell}(z),\\
\;_{r}e_{s}(A_{i}-nI)&=\prod_{k=1}^{n}(A_{i}-kI)\sum_{\ell=0}^{\infty}\prod_{k=1}^{n}\bigg{(}A_{i}+(\ell-k)I\bigg{)}^{-1}\Psi_{\ell}(z),\\
\;_{r}e_{s}(B_{j}+nI)&=\prod_{k=1}^{n}(B_{j}+(k-1)I)\sum_{\ell=0}^{\infty}\prod_{k=1}^{n}\bigg{(}B_{j}+(\ell+k-1)I\bigg{)}^{-1}\Psi_{\ell}(z),\\
\;_{r}e_{s}(B_{j}-nI)&=\prod_{k=1}^{n}\bigg{(}B_{j}-kI\bigg{)}^{-1}\sum_{\ell=0}^{\infty}\prod_{k=1}^{n}(B_{j}+(\ell-k)I)\Psi_{\ell}(z).\label{3.11}
\end{split}
\end{equation}
%---------------------------------------------------------------------------------------------------------------------------------------------------------------------------------------
%-----------------------------------------------------------------------------------------------------------------------------------------------------------------
\end{theorem}
%---------------------------------------------------------------------------------------------------------------------------------------------------------------------------------------
%---------------------------------------------------------------------------------------------------------------------------------------------------------------------------------------
\begin{proof}
%---------------------------------------------------------------------------------------------------------------------------------------------------------------------------------------
Using the relation $A_{i}(A_{i}+I)_{\ell}=(A_{i}+kI)(A_{i})_{\ell}$, we get the matrix contiguous function
%-----------------------------------------------------------------------------------------------------------------------------------------------------------------
\begin{equation}
\begin{split}
&\;_{r}e_{s}(A_{1}+)=\sum_{\ell=0}^{\infty}\frac{z^{\ell}}{n!}(A_{1}+I)_{\ell}(A_{2})_{\ell}\ldots(A_{r})_{\ell}[(B_{1})_{\ell}]^{-1}[(B_{2})_{\ell}]^{-1}\ldots[(B_{s})_{\ell}]^{-1}\gamma(\ell P+Q,x)\Gamma^{-1}(\ell P+Q)\\
&=\sum_{\ell=0}^{\infty}\frac{z^{\ell}}{n!}(A_{1}+\ell I)\bigg{(}A_{1}\bigg{)}^{-1}(A_{1})_{\ell}(A_{2})_{\ell}\ldots(A_{r})_{\ell}[(B_{1})_{\ell}]^{-1}[(B_{2})_{\ell}]^{-1}\ldots[(B_{s})_{\ell}]^{-1}\gamma(\ell P+Q,x)\Gamma^{-1}(\ell P+Q)\\
&=\sum_{\ell=0}^{\infty}(A_{1}+\ell I)\bigg{(}A_{1}\bigg{)}^{-1}\Psi_{\ell}(z).\label{3.12}
\end{split}
\end{equation}
%-----------------------------------------------------------------------------------------------------------------------------------------------------------------
Similarly, we get
%-----------------------------------------------------------------------------------------------------------------------------------------------------------------
\begin{equation}
\begin{split}
\;_{r}e_{s}(A_{i}+)&=\bigg{(}A_{i}\bigg{)}^{-1}\sum_{\ell=0}^{\infty}(A_{i}+\ell I)\Psi_{\ell}(z),\\
\;_{r}e_{s}(A_{i}-)&=(A_{i}-I)\sum_{\ell=0}^{\infty}\bigg{(}A_{i}+(\ell-1)I\bigg{)}^{-1}\Psi_{\ell}(z),\\
\;_{r}e_{s}(B_{j}+)&=(B_{j})\sum_{\ell=0}^{\infty}\bigg{(}B_{j}+\ell I\bigg{)}^{-1}\Psi_{\ell}(z),\\
\;_{r}e_{s}(B_{j}-)&=\bigg{(}B_{j}-I\bigg{)}^{-1}\sum_{\ell=0}^{\infty}(B_{j}+(k-1)I)\Psi_{\ell}(z).\label{3.13}
\end{split}
\end{equation}
%---------------------------------------------------------------------------------------------------------------------------------------------------------------------------------------
Repeating the previous steps, we get the equation (\ref{3.11}).
%---------------------------------------------------------------------------------------------------------------------------------------------------------------------------------------
\end{proof}
%---------------------------------------------------------------------------------------------------------------------------------------------------------------------------------------
%---------------------------------------------------------------------------------------------------------------------------------------------------------------------------------------
\begin{theorem} For the $\;_{r}e_{s}$ matrix function, the following recurrence relations holds true:
%---------------------------------------------------------------------------------------------------------------------------------------------------------------------------------------
%-----------------------------------------------------------------------------------------------------------------------------------------------------------------
\begin{equation}
\begin{split}
(A_{k}-A_{i})\;_{r}e_{s}=A_{k}\;_{r}e_{s}(A_{k}+)-A_{i}\;_{r}e_{s}(A_{i}+)\,;i\neq k;i,k=1,2,3,\ldots,r\label{3.14}
\end{split}
\end{equation}
%-----------------------------------------------------------------------------------------------------------------------------------------------------------------
%-----------------------------------------------------------------------------------------------------------------------------------------------------------------
\begin{equation}
\begin{split}
(B_{k}-B_{j})\;_{r}e_{s}=(B_{k}-I)\;_{r}e_{s}(B_{k}-)-(B_{j}-I)\;_{r}e_{s}(B_{j}-)\,;k\neq j;k,j=1,2,\ldots,s\label{3.15}
\end{split}
\end{equation}
%-----------------------------------------------------------------------------------------------------------------------------------------------------------------
and
%-----------------------------------------------------------------------------------------------------------------------------------------------------------------
\begin{equation}
\begin{split}
(A_{i}-B_{j}+I)\;_{r}e_{s}=A_{i}\;_{r}e_{s}(A_{i}+)-(B_{j}-I)\;_{r}e_{s}(B_{j}-)\,;i,k=1,2,3,\ldots,r,\,j=1,2,\ldots,s.\label{3.16}
\end{split}
\end{equation}
%---------------------------------------------------------------------------------------------------------------------------------------------------------------------------------------
\end{theorem}
%---------------------------------------------------------------------------------------------------------------------------------------------------------------------------------------
%---------------------------------------------------------------------------------------------------------------------------------------------------------------------------------------
\begin{proof}
%---------------------------------------------------------------------------------------------------------------------------------------------------------------------------------------
%-----------------------------------------------------------------------------------------------------------------------------------------------------------------
Using the differential operator $\theta=z\frac{d}{dz}$, we have
%-----------------------------------------------------------------------------------------------------------------------------------------------------------------
\begin{equation*}
\begin{split}
(\theta I+A_{i})\;_{r}e_{s}&=\sum_{\ell=0}^{\infty}(A_{i}+\ell I)\Psi_{\ell}(z).
\end{split}
\end{equation*}
%-----------------------------------------------------------------------------------------------------------------------------------------------------------------
Hence, with the aid of (\ref{3.13}), we get
%-----------------------------------------------------------------------------------------------------------------------------------------------------------------
\begin{equation}
\begin{split}
(\theta I+A_{i})\;_{r}e_{s}=A_{i}\;_{r}e_{s}(A_{i}+)\,;\,i=1,2,\ldots,r.\label{3.17}
\end{split}
\end{equation}
%-----------------------------------------------------------------------------------------------------------------------------------------------------------------
Similarly, we have
%-----------------------------------------------------------------------------------------------------------------------------------------------------------------
\begin{equation}
\begin{split}
(\theta I+B_{j}-I)\;_{r}e_{s}=(B_{j}-I)\;_{r}e_{s}(B_{j}-)\,;\,j=1,2,\ldots,s.\label{3.18}
\end{split}
\end{equation}
%-----------------------------------------------------------------------------------------------------------------------------------------------------------------
The elimination of $\theta\;_{r}e_{s}$ from (\ref{3.17}) and (\ref{3.18}) give the set of (\ref{3.14})-(\ref{3.16}).
%---------------------------------------------------------------------------------------------------------------------------------------------------------------------------------------
\end{proof}
%---------------------------------------------------------------------------------------------------------------------------------------------------------------------------------------
%---------------------------------------------------------------------------------------------------------------------------------------------------------------------------------------
\begin{theorem} For the $\;_{r}e_{s}$ matrix function, the following relation hold true:
%---------------------------------------------------------------------------------------------------------------------------------------------------------------------------------------
\begin{equation}
\begin{split}
D_{z}^{n}\;_{r}e_{s}=&\Pi_{i=1}^{r}(A_{i})_{n}\bigg{[}\Pi_{j=1}^{s}(B_{j})_{n}\bigg{]}^{-1}\;_{r}e_{s}(A_{1}+nI,\ldots,A_{r}+nI;B_{1}+nI,\ldots,B_{s}+nI;P,n P+Q,z).\label{3.19}
\end{split}
\end{equation}
%---------------------------------------------------------------------------------------------------------------------------------------------------------------------------------------
\end{theorem}
%---------------------------------------------------------------------------------------------------------------------------------------------------------------------------------------
%---------------------------------------------------------------------------------------------------------------------------------------------------------------------------------------
\begin{proof}
%---------------------------------------------------------------------------------------------------------------------------------------------------------------------------------------
It is well known that
%-----------------------------------------------------------------------------------------------------------------------------------------------------------------
\begin{equation}
\begin{split}
D_{z}\;_{r}e_{s}=&A_{1}\ldots\,A_{r}\bigg{(}B_{1}\bigg{)}^{-1}\ldots\bigg{(}B_{s}\bigg{)}^{-1}\;_{r}e_{s}(A_{1}+I,\\
&\ldots,A_{r}+I;B_{1}+I,\ldots,B_{s}+I;P,P+Q,z)\,;\\
=&\Pi_{i=1}^{r}A_{i}\bigg{(}\Pi_{j=1}^{s}B_{j}\bigg{)}^{-1}\;_{r}e_{s}(A_{1}+I,\\
&\ldots,A_{r}+I;B_{1}+I,\ldots,B_{s}+I;P,P+Q,z)\,;\, D_{z}=\frac{d}{dz}.\label{3.20}
\end{split}
\end{equation}
%-----------------------------------------------------------------------------------------------------------------------------------------------------------------
Repeating in the above equation, we get the equation (\ref{3.19}).
%-----------------------------------------------------------------------------------------------------------------------------------------------------------------
%---------------------------------------------------------------------------------------------------------------------------------------------------------------------------------------
\end{proof}
%---------------------------------------------------------------------------------------------------------------------------------------------------------------------------------------

%---------------------------------------------------------------------------------------------------------------------------------------------------------------------------------------------
\begin{theorem} Let $A$, $B$, $P$ and $Q$ be commutative matrices in $\Bbb{C}^{N\times N}$ satisfying the condition (\ref{2.1}). Then the following recursion formulas holds true for $\;_{r}e_{s}$
%---------------------------------------------------------------------------------------------------------------------------------------------------------------------------------------------
\begin{eqnarray}
\begin{split}
\bigg{(}\theta P+Q\bigg{)}\;_{r}e_{s}+e^{x}x^{Q}\;_{r}R_{s}(x^{P}z)=\bigg{(}\theta P+Q\bigg{)}\;_{r}e_{s}(Q+I).\label{4.12}
\end{split}
\end{eqnarray}
%---------------------------------------------------------------------------------------------------------------------------------------------------------------------------------------------
\end{theorem}
%---------------------------------------------------------------------------------------------------------------------------------------------------------------------------------------------
\begin{proof}
%---------------------------------------------------------------------------------------------------------------------------------------------------------------------------------------------
Starting with the right hand side, we have
%-----------------------------------------------------------------------------------------------------------------------------------------------------------------
\begin{equation*}
\begin{split}
&Q\;_{r}e_{s}(Q+I)+zP\frac{d}{dz}\;_{r}e_{s}(Q+I)\\
&=Q\;_{r}e_{s}(Q+I)+zP\frac{d}{dz}\bigg{[}\sum_{\ell=0}^{\infty}\frac{z^{\ell}}{\ell !}\prod_{i=1}^{r}(A_{i})_{\ell}\bigg{[}\prod_{j=1}^{s}(B_{j})_{\ell}\bigg{]}^{-1}\gamma(\ell P+Q+I,x)\Gamma^{-1}(\ell P+Q+I)\bigg{]}\\
&=Q\;_{r}e_{s}(Q+I)+zP\bigg{[}\sum_{\ell=0}^{\infty}\frac{\ell z^{\ell-1}}{\ell !}\prod_{i=1}^{r}(A_{i})_{\ell}\bigg{[}\prod_{j=1}^{s}(B_{j})_{\ell}\bigg{]}^{-1}\gamma(\ell P+Q+I,x)\Gamma^{-1}(\ell P+Q+I)\bigg{]}\\
&=Q\;_{r}e_{s}(Q+I)+\bigg{[}\sum_{\ell=0}^{\infty}\frac{(\ell P+Q-Q) z^{\ell}}{\ell !}\prod_{i=1}^{r}(A_{i})_{\ell}\bigg{[}\prod_{j=1}^{s}(B_{j})_{\ell}\bigg{]}^{-1}\gamma(\ell P+Q+I,x)\Gamma^{-1}(\ell P+Q+I)\bigg{]}
\end{split}
\end{equation*}
%---------------------------------------------------------------------------------------------------------------------------------------------------------------------------------------------
\begin{equation*}
\begin{split}
&=Q\;_{r}e_{s}(Q+I)+\sum_{\ell=0}^{\infty}\frac{(\ell P+Q) z^{\ell}}{\ell !}\prod_{i=1}^{r}(A_{i})_{\ell}\bigg{[}\prod_{j=1}^{s}(B_{j})_{\ell}\bigg{]}^{-1}\gamma(\ell P+Q+I,x)\Gamma^{-1}(\ell P+Q)(\ell P+Q)^{-1}\\
&-Q\sum_{\ell=0}^{\infty}\frac{ z^{\ell}}{\ell !}\prod_{i=1}^{r}(A_{i})_{\ell}\bigg{[}\prod_{j=1}^{s}(B_{j})_{\ell}\bigg{]}^{-1}\gamma(\ell P+Q+I,x)\Gamma^{-1}(\ell P+Q+I)\bigg{]}\\
&=\sum_{\ell=0}^{\infty}\frac{z^{\ell}}{\ell !}\prod_{i=1}^{r}(A_{i})_{\ell}\bigg{[}\prod_{j=1}^{s}(B_{j})_{\ell}\bigg{]}^{-1}\gamma(\ell P+Q+I,x)\Gamma^{-1}(\ell P+Q)\\
&=\sum_{\ell=0}^{\infty}\frac{z^{\ell}}{\ell !}\prod_{i=1}^{r}(A_{i})_{\ell}\bigg{[}\prod_{j=1}^{s}(B_{j})_{\ell}\bigg{]}^{-1}\bigg{[}(\ell P+Q)\gamma(\ell P+Q,x)+e^{x}x^{\ell P+Q}\bigg{]}\Gamma^{-1}(\ell P+Q)\\
&=\sum_{\ell=0}^{\infty}\frac{z^{\ell}}{\ell !}\prod_{i=1}^{r}(A_{i})_{\ell}\bigg{[}\prod_{j=1}^{s}(B_{j})_{\ell}\bigg{]}^{-1}(\ell P+Q)\gamma(\ell P+Q,x)\Gamma^{-1}(\ell P+Q)\\
&+\sum_{\ell=0}^{\infty}\frac{z^{\ell}}{\ell !}\prod_{i=1}^{r}(A_{i})_{\ell}\bigg{[}\prod_{j=1}^{s}(B_{j})_{\ell}\bigg{]}^{-1}e^{x}x^{\ell P+Q}\Gamma^{-1}(\ell P+Q)
\end{split}
\end{equation*}
%---------------------------------------------------------------------------------------------------------------------------------------------------------------------------------------------
\begin{equation*}
\begin{split}
&=P\sum_{\ell=0}^{\infty}\frac{\ell z^{\ell}}{\ell !}\prod_{i=1}^{r}(A_{i})_{\ell}\bigg{[}\prod_{j=1}^{s}(B_{j})_{\ell}\bigg{]}^{-1}\gamma(\ell P+Q,x)\Gamma^{-1}(\ell P+Q)\\
&+Q\sum_{\ell=0}^{\infty}\frac{z^{\ell}}{\ell !}\prod_{i=1}^{r}(A_{i})_{\ell}\bigg{[}\prod_{j=1}^{s}(B_{j})_{\ell}\bigg{]}^{-1}\gamma(\ell P+Q,x)\Gamma^{-1}(\ell P+Q)\\
&+e^{x}x^{Q}\sum_{\ell=0}^{\infty}\frac{x^{\ell P}z^{\ell}}{\ell !}\prod_{i=1}^{r}(A_{i})_{\ell}\bigg{[}\prod_{j=1}^{s}(B_{j})_{\ell}\bigg{]}^{-1}\Gamma^{-1}(\ell P+Q)\\
&=zP\frac{d}{dz}\;_{r}e_{s}+Q\;_{r}e_{s}+e^{x}x^{Q}\;_{r}R_{s}(x^{P}z).
\end{split}
\end{equation*}
%---------------------------------------------------------------------------------------------------------------------------------------------------------------------------------------------
\end{proof}
%---------------------------------------------------------------------------------------------------------------------------------------------------------------------------------------------
%---------------------------------------------------------------------------------------------------------------------------------------------------------------------------------------
\begin{theorem}
%---------------------------------------------------------------------------------------------------------------------------------------------------------------------------------------
The generalized incomplete exponential matrix functions satisfies the following integral representation:
%---------------------------------------------------------------------------------------------------------------------------------------------------------------------------------------
\begin{equation}
\begin{split}
{\displaystyle \;_{r}E_{s}\!\left[(x,P,Q;z)\left |{\begin{matrix} A_1,A_2,\dots,A_r\\ B_1,B_2,\dots,B_s\end{matrix}} \right.\right]} = \int_x^\infty t^{Q-1} \ e^{-t} \ {\displaystyle {}_rR_s\!\left[{\left.{\begin{matrix} A_1,A_2,\dots,A_r\\ B_1,B_2,\dots,B_s\end{matrix}} \right|P,Q;zt^P}\right]} dt\label{3.22}
\end{split}
\end{equation}
%---------------------------------------------------------------------------------------------------------------------------------------------------------------------------------------
%---------------------------------------------------------------------------------------------------------------------------------------------------------------------------------------
and
%---------------------------------------------------------------------------------------------------------------------------------------------------------------------------------------
\begin{equation}
\begin{split}
{\displaystyle \;_{r}e_{s}\!\left[(x,P,Q;z)\left |{\begin{matrix} A_1,A_2,\dots,A_r\\ B_1,B_2,\dots,B_s\end{matrix}} \right.\right]} = \int_{\infty}^{x} t^{Q-1} \ e^{-t} \ {\displaystyle {}_rR_s\!\left[{\left.{\begin{matrix} A_1,A_2,\dots,A_r\\ B_1,B_2,\dots,B_s\end{matrix}} \right|P,Q;zt^P}\right]} dt,\label{3.23}
\end{split}
\end{equation}
%---------------------------------------------------------------------------------------------------------------------------------------------------------------------------------------
%---------------------------------------------------------------------------------------------------------------------------------------------------------------------------------------
where $P,Q,A_i,B_j\in \mathbf{C}$ and $\Re(P)>0,\Re(Q)>0,\Re(A_i)>0,\Re(B_j)>0$, $\forall i=1,2, \dots ,r,\ \forall j=1,2,\dots,s$.
%---------------------------------------------------------------------------------------------------------------------------------------------------------------------------------------
\end{theorem}
%---------------------------------------------------------------------------------------------------------------------------------------------------------------------------------------
\begin{proof}
%---------------------------------------------------------------------------------------------------------------------------------------------------------------------------------------
On using the integral representation of incomplete gamma matrix function defined by (\ref{2.1}), we arrive at
%---------------------------------------------------------------------------------------------------------------------------------------------------------------------------------------
\begin{equation*}
\begin{split}
{\displaystyle \;_{r}E_{s}\!\left[(x,P,Q;z)\left |{\begin{matrix} A_1,\dots,A_r\\ B_1,\dots,B_s\end{matrix}} \right.\right]} =  \int_x^\infty  \ e^{-t} \left(\sum_{\ell=0}^{\infty}\frac{z^{\ell}}{\ell !}\prod_{i=1}^{r}(A_{i})_{\ell}\bigg{[}\prod_{j=1}^{s}(B_{j})_{\ell}\bigg{]}^{-1}\Gamma^{-1}(\ell P+Q) \right) t^{P \ell+Q-I}dt.
\end{split}
\end{equation*}
%---------------------------------------------------------------------------------------------------------------------------------------------------------------------------------------
Further simplification by reversing the order of summation and integration yields the r.h.s. of assertion (\ref{3.22}).
%---------------------------------------------------------------------------------------------------------------------------------------------------------------------------------------
\end{proof}
%---------------------------------------------------------------------------------------------------------------------------------------------------------------------------------------
\begin{corollary}
%---------------------------------------------------------------------------------------------------------------------------------------------------------------------------------------
By setting $P=I$ and $r=1,s=0$ i.e. $A_1=A$, (\ref{3.22}) and (\ref{3.23}) reduces to
%---------------------------------------------------------------------------------------------------------------------------------------------------------------------------------------
\begin{equation*}
\begin{split}
{\displaystyle _1e_0\!\left[(x, I, Q; z)\left |{\begin{matrix} A\\ - \end{matrix}} \right.\right]} = \Gamma^{-1}(Q) \int_{\infty}^{x} t^{Q-I} \ e^{-t} \ {\displaystyle {}_1F_1\!\left[{\left.{\begin{matrix} A\\ Q\end{matrix}} \ \right| zt}\right]} dt
\end{split}
\end{equation*}
%---------------------------------------------------------------------------------------------------------------------------------------------------------------------------------------
and
%---------------------------------------------------------------------------------------------------------------------------------------------------------------------------------------
\begin{equation*}
\begin{split}
{\displaystyle _1E_0\!\left[(x, I, Q; z)\left |{\begin{matrix} A\\ - \end{matrix}} \right.\right]} = \Gamma^{-1}(Q) \int_x^\infty t^{Q-I} \ e^{-t} \ {\displaystyle {}_1F_1\!\left[{\left.{\begin{matrix} A\\ Q\end{matrix}} \ \right| zt}\right]} dt.
\end{split}
\end{equation*}
%---------------------------------------------------------------------------------------------------------------------------------------------------------------------------------------
\end{corollary}
%---------------------------------------------------------------------------------------------------------------------------------------------------------------------------------------
\begin{corollary}
%---------------------------------------------------------------------------------------------------------------------------------------------------------------------------------------
For ${}\;_{r}R_{s}(P,Q;z)$ matrix function, the following integral representation holds true:
%---------------------------------------------------------------------------------------------------------------------------------------------------------------------------------------
\begin{equation*}
\begin{split}
{\displaystyle {}\;_{r}R_{s} \!\left[{\left.{\begin{matrix} A_1,A_2,\dots,A_r\\ B_1,B_2,\dots,B_s\end{matrix}} \right|P,Q;z}\right]} = \Gamma^{-1}(A_1) \int_0^\infty t^{A_1-1} e^{-t} \ {\displaystyle {}_{r-1}R_{s}\!\left[{\left.{\begin{matrix} A_2,\dots,A_r\\B_1,B_2,\dots,B_s\end{matrix}} \right|P,Q;zt}\right]} dt,
\end{split}
\end{equation*}
%---------------------------------------------------------------------------------------------------------------------------------------------------------------------------------------
where $r\leq s+1$, \ $\Re(A_1)>0$ and $\Re(P)>0,\Re(Q)>0$ and $\Re(A_i)>0,\Re(B_j)>0$, $\forall i=2, \dots ,r,\ \forall j=2,\dots,s$.
%---------------------------------------------------------------------------------------------------------------------------------------------------------------------------------------
\end{corollary}
%---------------------------------------------------------------------------------------------------------------------------------------------------------------------------------------
\begin{theorem}
%---------------------------------------------------------------------------------------------------------------------------------------------------------------------------------------
The generalized incomplete exponential $\;_{r}E_{s}(x,P,Q;z)$ matrix function have the following integral representation:
%---------------------------------------------------------------------------------------------------------------------------------------------------------------------------------------
\begin{equation}
\begin{split}
{\displaystyle \;_{r}e_{s}\!\left[(x,P,Q;z)\left |{\begin{matrix} A_1,A_2,\dots,A_r\\ B_1,B_2,\dots,B_s \end{matrix}} \right.\right]} &= \Gamma(B_{1}) \Gamma^{-1}(A_{1}) \Gamma^{-1}(B_{1}-A_{1}) \int_0^1 t^{A_1-1} (1-t)^{B_1-A_1-1} \ \\ & \times {\displaystyle _{r-1}e_{s-1}\!\left[(x,P,Q;zt)\left |{\begin{matrix} A_2,\dots,A_r\\ B_2,\dots,B_s\end{matrix}} \right.\right]} dt\label{3.24}
\end{split}
\end{equation}
%---------------------------------------------------------------------------------------------------------------------------------------------------------------------------------------
and
%---------------------------------------------------------------------------------------------------------------------------------------------------------------------------------------
\begin{equation}
\begin{split}
{\displaystyle \;_{r}E_{s}\!\left[(x,P,Q;z)\left |{\begin{matrix} A_1,A_2,\dots,A_r\\ B_1,B_2,\dots,B_s \end{matrix}} \right.\right]} &= \Gamma(B_{1}) \Gamma^{-1}(A_{1}) \Gamma^{-1}(B_{1}-A_{1}) \int_0^1 t^{A_1-1} (1-t)^{B_1-A_1-1} \ \\ & \times {\displaystyle _{r-1}E_{s-1}\!\left[(x,P,Q;zt)\left |{\begin{matrix} A_2,\dots,A_r\\ B_2,\dots,B_s\end{matrix}} \right.\right]} dt,\label{3.25}
\end{split}
\end{equation}
%---------------------------------------------------------------------------------------------------------------------------------------------------------------------------------------
where $\Re(A_1)>\Re(B_1)>0$ and $\Re(P)>0,\Re(Q)>0$ and $\Re(A_i)>0, \Re(B_j)>0$, $\forall i=2, \dots ,r,\ \forall j=2,\dots,s$.
%---------------------------------------------------------------------------------------------------------------------------------------------------------------------------------------
\end{theorem}
%---------------------------------------------------------------------------------------------------------------------------------------------------------------------------------------
\begin{proof}
%---------------------------------------------------------------------------------------------------------------------------------------------------------------------------------------
By adopting the following elementary integral definition of Beta matrix function:
%---------------------------------------------------------------------------------------------------------------------------------------------------------------------------------------
\begin{equation}
\begin{split}
(A_1)_{\ell}[(B_1)_{\ell}]^{-1}= \Gamma(B_{1}) \Gamma^{-1}(A_{1}) \Gamma^{-1}(B_{1}-A_{1}) \int_0^1 t^{A_1+(\ell-1)I} \ (1-t)^{B_1-A_1-I} \ dt,\label{3.26}
\end{split}
\end{equation}
%---------------------------------------------------------------------------------------------------------------------------------------------------------------------------------------
in the l.h.s. of (\ref{3.24}), this yields
%---------------------------------------------------------------------------------------------------------------------------------------------------------------------------------------
\begin{equation*}
\begin{split}
{\displaystyle \;_{r}e_{s}\!\left[(x,P,Q;z)\left |{\begin{matrix} A_1,A_2,\dots,A_r\\ B_1,B_2,\dots,B_s \end{matrix}} \right.\right]}
&=  \Gamma(B_{1}) \Gamma^{-1}(A_{1}) \Gamma^{-1}(B_{1}-A_{1}) \sum_{\ell=0}^{\infty}\frac{z^\ell}{\ell!} \prod_{i=1}^{r}(A_{i})_{\ell}\bigg{[}\prod_{j=1}^{s}(B_{j})_{\ell}\bigg{]}^{-1}\\ & \times \gamma(P\ell+Q,x)\Gamma^{-1}(P\ell+Q)\int_0^1 t^{A_1+(\ell-1)I} \ (1-t)^{B_1-A_1-1} \ dt.
\end{split}
\end{equation*}
%---------------------------------------------------------------------------------------------------------------------------------------------------------------------------------------
Further simplification by reversing the order of summation and integration yields the r.h.s. of assertion (\ref{3.24}).
%---------------------------------------------------------------------------------------------------------------------------------------------------------------------------------------
\end{proof}
%---------------------------------------------------------------------------------------------------------------------------------------------------------------------------------------
\begin{corollary}[\cite{ds}]
%---------------------------------------------------------------------------------------------------------------------------------------------------------------------------------------
For the ${}\;_{r}R_{s}(P,Q;z)$ matrix function, we have the following integral representation:
%---------------------------------------------------------------------------------------------------------------------------------------------------------------------------------------
%---------------------------------------------------------------------------------------------------------------------------------------------------------------------------------------
\begin{equation*}
\begin{split}
{\displaystyle {}\;_{r}R_{s} \!\left[{\left.{\begin{matrix} A_1,A_2,\dots,A_r\\ B_1,B_2,\dots,B_s\end{matrix}} \right|P,Q;z}\right]} &= \Gamma(B_{1}) \Gamma^{-1}(A_{1}) \Gamma^{-1}(B_{1}-A_{1}) \int_0^1 t^{A_1-1} (1-t)^{B_1-A_1-I} \\ & \times {\displaystyle {}_{r-1}R_{s-1}\!\left[{\left.{\begin{matrix} A_2,\dots,A_r\\ B_2,\dots,B_s\end{matrix}} \right|P,Q;zt}\right]} dt,
\end{split}
\end{equation*}
%---------------------------------------------------------------------------------------------------------------------------------------------------------------------------------------
where $r\leq s+1$, \ $\Re(B_1)>\Re(A_1)>0$ and $\Re(P)>0,\Re(Q)>0$ and $\Re(A_i)>0,\Re(B_j)>0$, $\forall i=2, \dots ,r,\ \forall j=2,\dots,s$.
%---------------------------------------------------------------------------------------------------------------------------------------------------------------------------------------
\end{corollary}
%---------------------------------------------------------------------------------------------------------------------------------------------------------------------------------------
\begin{theorem}
%---------------------------------------------------------------------------------------------------------------------------------------------------------------------------------------
The generalized incomplete exponential $\;_{r}E_{s}(x,P,Q;z)$ matrix function have the following derivative formula:
%---------------------------------------------------------------------------------------------------------------------------------------------------------------------------------------
%---------------------------------------------------------------------------------------------------------------------------------------------------------------------------------------
\begin{equation}
\begin{split}
\frac{d^n}{dz^n} \left\lbrace {\displaystyle \;_{r}E_{s}\!\left[(x,P,Q;z)\left |{\begin{matrix} A_1,A_2,\dots,A_r\\ B_1,B_2,\dots,B_s\end{matrix}} \right.\right]} \right\rbrace &= \prod_{i=1}^{r}(A_{i})_{n}\bigg{[}\prod_{j=1}^{s}(B_{j})_{n}\bigg{]}^{-1} \\ & \times {\displaystyle \;_{r}E_{s}\!\left[(x,P,nP+Q;z)\left |{\begin{matrix} A_1+nI,\dots,A_r+nI\\ B_1+nI,\dots,B_s+nI\end{matrix}} \right.\right]} \label{3.27}
\end{split}
\end{equation}
%---------------------------------------------------------------------------------------------------------------------------------------------------------------------------------------
and
%---------------------------------------------------------------------------------------------------------------------------------------------------------------------------------------
\begin{equation}
\begin{split}
\frac{d^n}{dz^n} \left\lbrace {\displaystyle \;_{r}e_{s}\!\left[(x,P,Q;z)\left |{\begin{matrix} A_1,A_2,\dots,A_r\\ B_1,B_2,\dots,B_s\end{matrix}} \right.\right]} \right\rbrace &= \prod_{i=1}^{r}(A_{i})_{n}\bigg{[}\prod_{j=1}^{s}(B_{j})_{n}\bigg{]}^{-1} \\ & \times {\displaystyle \;_{r}e_{s}\!\left[(x,P,nP+Q;z)\left |{\begin{matrix} A_1+nI,\dots,A_r+nI\\ B_1+nI,\dots,B_s+nI\end{matrix}} \right.\right]}, \label{3.27}
\end{split}
\end{equation}
%---------------------------------------------------------------------------------------------------------------------------------------------------------------------------------------
%---------------------------------------------------------------------------------------------------------------------------------------------------------------------------------------
where $P,Q,A_i,B_j\in \mathbf{C}$ and $\Re(P)>0,\Re(Q)>0$ and $\Re(A_i)>0, \Re(B_j)>0$, $\forall i=1,2, \dots ,r,\ \forall j=1,2,\dots,s$.
%---------------------------------------------------------------------------------------------------------------------------------------------------------------------------------------
\end{theorem}
%---------------------------------------------------------------------------------------------------------------------------------------------------------------------------------------
\begin{proof}
%---------------------------------------------------------------------------------------------------------------------------------------------------------------------------------------
Differentiating (\ref{3.6}) with respect to $z$ and replacing $\ell$ by $\ell+1$, we arrive at
%---------------------------------------------------------------------------------------------------------------------------------------------------------------------------------------
%---------------------------------------------------------------------------------------------------------------------------------------------------------------------------------------
\begin{equation*}
\begin{split}
\frac{d}{dz} \left\lbrace {\displaystyle \;_{r}E_{s}\!\left[(x,P,Q;z)\left |{\begin{matrix} A_1,A_2,\dots,A_r\\ B_1,B_2,\dots,B_s\end{matrix}} \right.\right]} \right\rbrace \\
=\sum_{\ell=0}^{\infty}\frac{z^{\ell}}{\ell !}\prod_{i=1}^{r}(A_{i})_{\ell+1}\bigg{[}\prod_{j=1}^{s}(B_{j})_{\ell+1}\bigg{]}^{-1}\Gamma(\ell P+P+Q,x)\Gamma^{-1}(\ell P+P+Q).
\end{split}
\end{equation*}
%---------------------------------------------------------------------------------------------------------------------------------------------------------------------------------------
Using the relation $(A)_{\ell+1}=A(A+I)_\ell$, we get
%---------------------------------------------------------------------------------------------------------------------------------------------------------------------------------------
\begin{equation*}
\begin{split}
\frac{d}{dz} \left\lbrace {\displaystyle \;_{r}E_{s}\!\left[(x,P,Q;z)\left |{\begin{matrix} A_1,A_2,\dots,A_r\\ B_1,B_2,\dots,B_s\end{matrix}} \right.\right]} \right\rbrace = \prod_{i=1}^{r}(A_{i})\bigg{[}\prod_{j=1}^{s}(B_{j})\bigg{]}^{-1}\times {\displaystyle \;_{r}E_{s}\!\left[(x,P,P+Q;z)\left |{\begin{matrix} A_1+I,\dots,A_r+I \\ B_1+I,\dots,B_s+I\end{matrix}} \right.\right]}.
\end{split}
\end{equation*}
%---------------------------------------------------------------------------------------------------------------------------------------------------------------------------------------
%---------------------------------------------------------------------------------------------------------------------------------------------------------------------------------------
By repeating above procedure $n$-times yields the r.h.s. of assertion (\ref{3.27}).
%---------------------------------------------------------------------------------------------------------------------------------------------------------------------------------------
\end{proof}
%---------------------------------------------------------------------------------------------------------------------------------------------------------------------------------------------
\begin{corollary}
%---------------------------------------------------------------------------------------------------------------------------------------------------------------------------------------
For the ${}\;_{r}R_{s}(P,Q;z)$ matrix function, we have the following derivative formula:
%---------------------------------------------------------------------------------------------------------------------------------------------------------------------------------------
%---------------------------------------------------------------------------------------------------------------------------------------------------------------------------------------------
\begin{equation*}
\begin{split}
\frac{d^n}{dz^n} \left\lbrace {\displaystyle {}_rR_s\!\left[{\left.{\begin{matrix} A_1,A_2,\dots,a_p\\ B_1,B_2,\dots,B_s\end{matrix}} \right|P,Q;z}\right]} \right\rbrace &= \prod_{i=1}^{r}(A_{i})_{\ell}\bigg{[}\prod_{j=1}^{s}(B_{j})_{\ell}\bigg{]}^{-1} \\ & \times {\displaystyle {}_rR_s\!\left[{\left.{\begin{matrix} A_1+nI,A_2+nI,\dots,A_r+nI\\ B_1+nI,B_2+nI,\dots,B_s+nI\end{matrix}} \right|P,nP+Q;z}\right]},
\end{split}
\end{equation*}
%---------------------------------------------------------------------------------------------------------------------------------------------------------------------------------------------
where $P,Q,A_i,B_j\in \mathbf{C}$ and $\Re(P)>0,\Re(Q)>0$ and $\Re(A_i)>0, \Re(B_j)>0$, $\forall i=1,2, \dots ,r,\ \forall j=1,2,\dots,s$.
%---------------------------------------------------------------------------------------------------------------------------------------------------------------------------------------
\end{corollary}
%---------------------------------------------------------------------------------------------------------------------------------------------------------------------------------------
\begin{theorem}
%---------------------------------------------------------------------------------------------------------------------------------------------------------------------------------------
For the generalized incomplete exponential matrix functions, the following partial derivatives holds true:
%---------------------------------------------------------------------------------------------------------------------------------------------------------------------------------------
%---------------------------------------------------------------------------------------------------------------------------------------------------------------------------------------------
\begin{equation}
\begin{split}
\frac{\partial}{\partial z} \left\lbrace {\displaystyle \;_{r}E_{s}\!\left[(x,P,Q;z)\left |{\begin{matrix} A_1,A_2,\dots,A_r\\ B_1,B_2,\dots,B_s\end{matrix}} \right.\right]} \right\rbrace &= A_1\dots A_r[B_1]^{-1}\dots [B_s]^{-1} \\ & \times {\displaystyle \;_{r}E_{s}\!\left[(x,P,P+Q;z)\left |{\begin{matrix} A_1+I,\dots,A_r+I\\ B_1+1,\dots,B_s+I\end{matrix}} \right.\right]},\label{3.29}
\end{split}
\end{equation}
%---------------------------------------------------------------------------------------------------------------------------------------------------------------------------------------------
%---------------------------------------------------------------------------------------------------------------------------------------------------------------------------------------------
\begin{equation}
\begin{split}
\frac{\partial}{\partial z} \left\lbrace {\displaystyle \;_{r}e_{s}\!\left[(x,P,Q;z)\left |{\begin{matrix} A_1,A_2,\dots,A_r\\ B_1,B_2,\dots,B_s\end{matrix}} \right.\right]} \right\rbrace &= A_1\dots A_r[B_1]^{-1}\dots [B_s]^{-1} \\ & \times {\displaystyle \;_{r}e_{s}\!\left[(x,P,P+Q;z)\left |{\begin{matrix} A_1+I,\dots,A_r+I\\ B_1+1,\dots,B_s+I\end{matrix}} \right.\right]},\label{3.30}
\end{split}
\end{equation}
%---------------------------------------------------------------------------------------------------------------------------------------------------------------------------------------------
%---------------------------------------------------------------------------------------------------------------------------------------------------------------------------------------------
\begin{equation}
\begin{split}
\frac{\partial}{\partial x} \left\lbrace {\displaystyle \;_{r}e_{s}\!\left[(x,P,Q;z)\left |{\begin{matrix} A_1,A_2,\dots,A_r\\ B_1,B_2,\dots,B_s\end{matrix}} \right.\right]} \right\rbrace &=  e^{-x} x^{Q-I} {\displaystyle {}_rR_s\!\left[{\left.{\begin{matrix} A_1,A_2,\dots,A_r\\ B_1,B_2,\dots,B_s\end{matrix}} \right|P,Q;zx^P}\right]} ,\label{3.31}
\end{split}
\end{equation}
%---------------------------------------------------------------------------------------------------------------------------------------------------------------------------------------------
and
%---------------------------------------------------------------------------------------------------------------------------------------------------------------------------------------------
%---------------------------------------------------------------------------------------------------------------------------------------------------------------------------------------------
\begin{equation}
\begin{split}
\frac{\partial}{\partial x} \left\lbrace {\displaystyle \;_{r}E_{s}\!\left[(x,P,Q;z)\left |{\begin{matrix} A_1,A_2,\dots,A_r\\ B_1,B_2,\dots,B_s\end{matrix}} \right.\right]} \right\rbrace &= - e^{-x} x^{Q-I} {\displaystyle {}_rR_s\!\left[{\left.{\begin{matrix} A_1,A_2,\dots,A_r\\ B_1,B_2,\dots,B_s\end{matrix}} \right|P,Q;zx^P}\right]} ,\label{3.32}
\end{split}
\end{equation}
%---------------------------------------------------------------------------------------------------------------------------------------------------------------------------------------------
where $P,Q,A_i,B_j\in \mathbf{C}$ and $\Re(P)>0,\Re(Q)>0$ and $\Re(A_i)>0,\Re(B_j)>0$, $\forall i=1,2, \dots ,r,\ \forall j=1,2,\dots,s$.
%---------------------------------------------------------------------------------------------------------------------------------------------------------------------------------------
\end{theorem}
%---------------------------------------------------------------------------------------------------------------------------------------------------------------------------------------
\begin{proof}
%---------------------------------------------------------------------------------------------------------------------------------------------------------------------------------------
Differentiating partially (\ref{3.7}) with respect to $z$ and treating $x$ as a constant, we have
%---------------------------------------------------------------------------------------------------------------------------------------------------------------------------------------
%---------------------------------------------------------------------------------------------------------------------------------------------------------------------------------------------
\begin{equation*}
\begin{split}
&\frac{\partial}{\partial z} \left\lbrace {\displaystyle \;_{r}E_{s}\!\left[(x,P,Q;z)\left |{\begin{matrix} A_1,A_2,\dots,A_r\\ B_1,B_2,\dots,B_s\end{matrix}} \right.\right]} \right\rbrace \\
 &= \frac{\partial}{\partial z} \left\lbrace \sum_{\ell=0}^{\infty} \Gamma(P \ell+Q,x)\Gamma^{-1}(P \ell+Q)(A_1)_\ell(A_2)_\ell\dots(A_r)_\ell[(B_1)_\ell]^{-1}[(B_2)_\ell]^{-1}\dots[(B_s)_\ell]^{-1} \frac{z^\ell}{\ell !} \right\rbrace ,\\
&= \sum_{\ell=1}^{\infty} \Gamma(P \ell+Q,x)\Gamma^{-1}(P \ell+Q)(A_1)_\ell(A_2)_\ell\dots(A_r)_\ell[(B_1)_\ell]^{-1}[(B_2)_\ell]^{-1}\dots[(B_s)_\ell]^{-1} \frac{z^{\ell-1}}{(\ell-1)!}.
\end{split}
\end{equation*}
%---------------------------------------------------------------------------------------------------------------------------------------------------------------------------------------------
%---------------------------------------------------------------------------------------------------------------------------------------------------------------------------------------
By replacing $\ell$ by $\ell+1$, this leads to proof of (\ref{3.29}). For the proof of (\ref{3.31}), we differentiate partially first integral representation (\ref{3.7}) with respect to $x$ and treating $z$ as a constant.
%---------------------------------------------------------------------------------------------------------------------------------------------------------------------------------------
\end{proof}
%---------------------------------------------------------------------------------------------------------------------------------------------------------------------------------------
\begin{theorem}\label{thm3.9}
%---------------------------------------------------------------------------------------------------------------------------------------------------------------------------------------
For the generalized incomplete exponential $\;_{r}E_{s}(x,P,Q;z)$ matrix function, the following addition formula for addition of two argument is valid:
%---------------------------------------------------------------------------------------------------------------------------------------------------------------------------------------
%---------------------------------------------------------------------------------------------------------------------------------------------------------------------------------------------
\begin{equation}
\begin{split}
{\displaystyle \;_{r}E_{s}\!\left[(x,P,Q;y+z)\left |{\begin{matrix} A_1,A_2,\dots,A_r\\ B_1,B_2,\dots,B_s\end{matrix}} \right.\right]}  &=  \Gamma(B_1)\dots \Gamma(B_s)\Gamma^{-1}(A_1)\dots \Gamma^{-1}(A_r)\sum_{\ell=0}^\infty \Gamma(A_1+\ell I)\dots \Gamma(A_r+\ell I) \\
&\Gamma^{-1}(B_1+\ell I)\dots \Gamma^{-1}(B_s+\ell I) \frac{z^\ell }{\ell !}\\
&{\displaystyle \;_{r}E_{s}\!\left[(x,P,P \ell+Q;y)\left |{\begin{matrix} A_1+\ell I,\dots,A_r+\ell I\\ B_1+\ell I,\dots,B_s+\ell I\end{matrix}} \right.\right]}.\label{3.34}
\end{split}
\end{equation}
%---------------------------------------------------------------------------------------------------------------------------------------------------------------------------------------------
%---------------------------------------------------------------------------------------------------------------------------------------------------------------------------------------
\end{theorem}
%---------------------------------------------------------------------------------------------------------------------------------------------------------------------------------------
\begin{proof}
%---------------------------------------------------------------------------------------------------------------------------------------------------------------------------------------
Proof of above theorem pursue from the derivative formula of $\;_{r}E_{s}(x,P,Q;z)$ matrix function:
%---------------------------------------------------------------------------------------------------------------------------------------------------------------------------------------
%---------------------------------------------------------------------------------------------------------------------------------------------------------------------------------------------
\begin{equation*}
\begin{split}
\frac{d^n}{dz^n} \left\lbrace {\displaystyle \;_{r}E_{s}\!\left[(x,P,Q;z)\left |{\begin{matrix} A_1,A_2,\dots,A_r\\ B_1,B_2,\dots,B_s\end{matrix}} \right.\right]} \right\rbrace &= \prod_{i=1}^{r}(A_{i})_{n}\bigg{[}\prod_{j=1}^{s}(B_{j})_{n}\bigg{]}^{-1} \\ & \times {\displaystyle \;_{r}E_{s}\!\left[(x,P,n P+Q;z)\left | \ {\begin{matrix} A_1+nI,\dots,A_r+nI\\ B_1+nI,\dots,B_s+nI\end{matrix}} \right.\right]}.
\end{split}
\end{equation*}
%---------------------------------------------------------------------------------------------------------------------------------------------------------------------------------------------
\end{proof}
%---------------------------------------------------------------------------------------------------------------------------------------------------------------------------------------------
\begin{theorem} \label{thm3.10}
%---------------------------------------------------------------------------------------------------------------------------------------------------------------------------------------
For the generalized incomplete exponential $\;_{r}E_{s}(x,P,Q;z)$ matrix function, the following multiplication formula for multiplication of two argument is valid:
%---------------------------------------------------------------------------------------------------------------------------------------------------------------------------------------
%---------------------------------------------------------------------------------------------------------------------------------------------------------------------------------------------
\begin{equation}
\begin{split}
{\displaystyle \;_{r}E_{s}\!\left[(x,P,Q;yz)\left |{\begin{matrix} A_1,A_2,\dots,A_r\\ B_1,B_2,\dots,B_s\end{matrix}} \right.\right]}  &= \Gamma(B_1)\dots \Gamma(B_s)\Gamma^{-1}(A_1)\dots \Gamma^{-1}(A_r) \sum_{n=0}^\infty \Gamma(A_1+nI)\dots \Gamma(A_r+nI) \\
&\Gamma^{-1}(B_1+nI)\dots \Gamma^{-1}(B_s+nI)\frac{y^n(z-1)^n}{n!} \\
&{\displaystyle \;_{r}E_{s}\!\left[(x,P,P n+Q;y)\left |{\begin{matrix} A_1+nI,\dots,A_r+nI \\ B_1+nI,\dots,B_s+nI\end{matrix}} \right.\right]}.\label{3.35}
\end{split}
\end{equation}
%---------------------------------------------------------------------------------------------------------------------------------------------------------------------------------------------
%---------------------------------------------------------------------------------------------------------------------------------------------------------------------------------------
\end{theorem}
%---------------------------------------------------------------------------------------------------------------------------------------------------------------------------------------
\begin{proof}
%---------------------------------------------------------------------------------------------------------------------------------------------------------------------------------------
Proof of above theorem is similar to Theorem \ref{thm3.9}.
%---------------------------------------------------------------------------------------------------------------------------------------------------------------------------------------
\end{proof}
%---------------------------------------------------------------------------------------------------------------------------------------------------------------------------------------
\begin{theorem}
%---------------------------------------------------------------------------------------------------------------------------------------------------------------------------------------
For the generalized incomplete exponential $\;_{r}E_{s}(x,P,Q;z)$ matrix function, we have the following integral representation:
%---------------------------------------------------------------------------------------------------------------------------------------------------------------------------------------
%---------------------------------------------------------------------------------------------------------------------------------------------------------------------------------------------
\begin{equation}
\begin{split}
&\int_0^t u^{P-1} (t-u)^{Q-I} {\displaystyle \;_{r}E_{s}\!\left[(x,P,Q; z u^k)\left |{\begin{matrix} A_1,A_2,\dots,A_r\\ B_1,B_2,\dots,B_s\end{matrix}} \right.\right]} du = \mathbf{B}(P,Q) \ t^{P+Q-I} \\ & \;\;\;\;\;\;\;\;\;\;\;\;\;\;\;\;\;\;\;\;\;\;\;\;\;\;\;\;\;\;\;\;\;\;\;\;\;\;\; \times {\displaystyle _{r+k}E_{s+k}\!\left[(x,P,Q; z t^k)\left |{\begin{matrix} A_1,A_2,\dots,A_r, \Delta \left( {k, P} \right) \\ B_1,B_2,\dots,B_s,\Delta \left( {k, P+Q} \right)\end{matrix}} \right.\right]}, \label{4.25}
\end{split}
\end{equation}
%---------------------------------------------------------------------------------------------------------------------------------------------------------------------------------------------

where $\Delta \left( {k, P} \right)$ represents the sequence of $k$ parameters i.e.
%---------------------------------------------------------------------------------------------------------------------------------------------------------------------------------------------
\begin{equation*}
\begin{split}
\frac{P}{k},\frac{P  + I}{k },\frac{P  + 2I}{k}, \ldots ,\frac{{ P + (k - 1)I}}{k },
\end{split}
\end{equation*}
%---------------------------------------------------------------------------------------------------------------------------------------------------------------------------------------------
and $\Re(P)>0,\Re(Q)>0$ and $\Re(A_i)>0,\Re(B_j)>0$, $\forall i=1,2, \dots ,r,\ \forall j=1,2,\dots,s$.
%---------------------------------------------------------------------------------------------------------------------------------------------------------------------------------------------
\end{theorem}
%---------------------------------------------------------------------------------------------------------------------------------------------------------------------------------------------
\begin{proof}
%---------------------------------------------------------------------------------------------------------------------------------------------------------------------------------------------
Let $\mathcal{R}$ be the l.h.s. of (\ref{4.25}). Then, using (\ref{3.7}), this gives
%---------------------------------------------------------------------------------------------------------------------------------------------------------------------------------------------
\begin{equation*}
\begin{split}
\mathcal{R} =& \int_0^t u^{P-1} (t-u)^{P-I}  \sum_{\ell=0}^{\infty} \Gamma(P \ell+Q,x)\Gamma^{-1}(P\ell+Q)\\
&(A_1)_\ell(A_2)_\ell\dots(A_r)_\ell[(B_1)_\ell]^{-1}[(B_2)_\ell]^{-1}\dots[(B_s)_\ell]^{-1} \frac{(z u^k)^\ell}{\ell !}  du.
\end{split}
\end{equation*}
%---------------------------------------------------------------------------------------------------------------------------------------------------------------------------------------------

%---------------------------------------------------------------------------------------------------------------------------------------------------------------------------------------
Substituting $u=tx$, we arrive at
%---------------------------------------------------------------------------------------------------------------------------------------------------------------------------------------
%---------------------------------------------------------------------------------------------------------------------------------------------------------------------------------------------
\begin{equation*}
\begin{split}
\mathcal{R} &= t^{P+Q-I}\sum_{\ell=0}^{\infty} \int_0^1 x^{P+(k\ell-1)I} (1-x)^{P-I}dx   \Gamma(P \ell+Q,x)\Gamma^{-1}(P\ell+Q)\\
&(A_1)_\ell(A_2)_\ell\dots(A_r)_\ell[(B_1)_\ell]^{-1}[(B_2)_\ell]^{-1}\dots[(B_s)_\ell]^{-1} \frac{(z t^k)^\ell}{\ell !}\\
&= t^{P+Q-I} \sum_{\ell=0}^{\infty} \Gamma(P \ell+Q,x)\Gamma^{-1}(P\ell+Q)\\
&(A_1)_\ell(A_2)_\ell\dots(A_r)_\ell[(B_1)_\ell]^{-1}[(B_2)_\ell]^{-1}\dots[(B_s)_\ell]^{-1} \frac{(z t^k)^\ell}{\ell !} \; \Gamma(P+k \ell)\Gamma(Q)\Gamma^{-1}(P+Q+k\ell).
\end{split}
\end{equation*}
%---------------------------------------------------------------------------------------------------------------------------------------------------------------------------------------------
%---------------------------------------------------------------------------------------------------------------------------------------------------------------------------------------
Now using the property of Pochhammer symbol defined in (\ref{2.11}), this leads to the right hand side of (\ref{4.25}).
%---------------------------------------------------------------------------------------------------------------------------------------------------------------------------------------
\end{proof}
%---------------------------------------------------------------------------------------------------------------------------------------------------------------------------------------
\begin{theorem}
%---------------------------------------------------------------------------------------------------------------------------------------------------------------------------------------
For the generalized incomplete exponential $\;_{r}E_{s}(x,P,Q;z)$ matrix function, the following integral representation holds true:
%---------------------------------------------------------------------------------------------------------------------------------------------------------------------------------------
%---------------------------------------------------------------------------------------------------------------------------------------------------------------------------------------------
\begin{equation}
\begin{split}
& \Gamma^{-1}(Q)\Gamma^{-1}(C)\Gamma(Q+C)\int_t^x (x-u)^{C -1} (u-t)^{Q-I} {\displaystyle \;_{r}E_{s}\!\left[(x,P,Q; z (u-t)^k)\left |{\begin{matrix} A_1,A_2,\dots,A_r\\ B_1,B_2,\dots,B_s\end{matrix}} \right.\right]} du    \\ & \;\;\;\;\;\;\;\;\;\;\;\;\; =(x-t)^{C+Q-I} {\displaystyle _{r+k}E_{s+k}\!\left[(x,P,Q; z (x-t)^k)\left |{\begin{matrix} A_1,A_2,\dots,A_r,\Delta \left( {k, Q} \right)\\ B_1,B_2,\dots,B_s,\Delta \left( {k, Q+C} \right)\end{matrix}} \right.\right]}, \label{3.37}
\end{split}
\end{equation}
%---------------------------------------------------------------------------------------------------------------------------------------------------------------------------------------------
%---------------------------------------------------------------------------------------------------------------------------------------------------------------------------------------
where $\Delta \left( {k, Q} \right)$ represents the sequence of $k$ parameters i.e.
%---------------------------------------------------------------------------------------------------------------------------------------------------------------------------------------
%---------------------------------------------------------------------------------------------------------------------------------------------------------------------------------------------
\begin{equation*}
\begin{split}
\frac{Q}{k},\frac{Q  + I}{k },\frac{Q  + 2I}{k }, \ldots ,\frac{{ Q +( k - 1)I}}{k },
\end{split}
\end{equation*}
%-----------------------------------------------------------------------------------------------------------------------------------------------------------------------------------
%---------------------------------------------------------------------------------------------------------------------------------------------------------------------------------------
and $\Re(P)>0,\Re(Q)>0$ and $\Re(A_i)>0, \Re(B_j)>0$, $\forall i=1,2, \dots ,r,\ \forall j=1,2,\dots,s$.
%---------------------------------------------------------------------------------------------------------------------------------------------------------------------------------------
\end{theorem}
%---------------------------------------------------------------------------------------------------------------------------------------------------------------------------------------
\begin{proof}
%---------------------------------------------------------------------------------------------------------------------------------------------------------------------------------------
Let $\mathbf{R}$ be the l.h.s. of (\ref{3.37}), and using (\ref{3.7}), we get
%---------------------------------------------------------------------------------------------------------------------------------------------------------------------------------------
%---------------------------------------------------------------------------------------------------------------------------------------------------------------------------------------------
\begin{equation*}
\begin{split}
\mathbf{R}=& \Gamma^{-1}(Q)\Gamma^{-1}(C)\Gamma(Q+C) \int_t^x (x-u)^{C -1} (u-t)^{Q-I}  \sum_{n=0}^{\infty} \Gamma(P n+Q,x)\Gamma^{-1}(Pn+Q)\\
&(A_1)_n(A_2)_n\dots(A_r)_n[(B_1)_n]^{-1}[(B_2)_n]^{-1}\dots[(B_s)_n]^{-1} \frac{(z (u-t)^k)^n}{n!}  du.
\end{split}
\end{equation*}
%-----------------------------------------------------------------------------------------------------------------------------------------------------------------------------------
%---------------------------------------------------------------------------------------------------------------------------------------------------------------------------------------
Substituting $y=\frac{u-t}{x-t}$, we arrive at
%---------------------------------------------------------------------------------------------------------------------------------------------------------------------------------------
%---------------------------------------------------------------------------------------------------------------------------------------------------------------------------------------------
\begin{equation*}
\begin{split}
\mathbf{R}&= (x-t)^{C+Q-I}\Gamma^{-1}(Q)\Gamma^{-1}(C)\Gamma(Q+C) \int_0^1 y^{kn I+Q-I} \ (1-y)^{C-1}  \sum_{n=0}^{\infty} \Gamma(P n+Q,x)\Gamma^{-1}(Pn+Q)\\
&(A_1)_n(A_2)_n\dots(A_r)_n[(B_1)_n]^{-1}[(B_2)_n]^{-1}\dots[(B_s)_n]^{-1}\frac{(z(x-t)^k)^n}{n!} \ dv, \\
&= (x-t)^{C+Q-I}\Gamma^{-1}(Q)\Gamma^{-1}(C)\Gamma(Q+C) \sum_{n=0}^{\infty} \Gamma(P n+Q,x)\Gamma^{-1}(Pn+Q)(A_1)_n(A_2)_n\dots(A_r)_n\\
&[(B_1)_n]^{-1}[(B_2)_n]^{-1}\dots[(B_s)_n]^{-1} \frac{(z(x-t)^k)^n}{n!}  \; \Gamma(kn+Q)\Gamma(C)\Gamma^{-1}(kn+Q+C )\\
&= (x-t)^{C+Q-I}\sum_{n=0}^{\infty} \Gamma(P n+Q,x)\Gamma^{-1}(Pn+Q)(A_1)_n(A_2)_n\dots(A_r)_n\\
&[(B_1)_n]^{-1}[(B_2)_n]^{-1}\dots[(B_s)_n]^{-1} \frac{(z(x-t)^k)^n}{n!}   (Q)_{kn}[(Q+C)_{kn}]^{-1}.
\end{split}
\end{equation*}
%-----------------------------------------------------------------------------------------------------------------------------------------------------------------------------------
Now using the property of Pochhammer symbol defined in (\ref{2.11}), this leads to the right hand side of (\ref{3.37}).
%---------------------------------------------------------------------------------------------------------------------------------------------------------------------------------------
\end{proof}
%---------------------------------------------------------------------------------------------------------------------------------------------------------------------------------------
%-----------------------------------------------------------------------------------------------------------------------------------------------------------------------------------
\begin{theorem}
%-----------------------------------------------------------------------------------------------------------------------------------------------------------------------------------
Suppose that $A_{i},B_{j}$ and $B_{j}-A_{i}$($1\leq i\leq r$, $B_{j}$\,;\,$1\leq j\leq s$) are positive stable matrices in $\Bbb{C}^{N \times N}$ such that $B_{j}+\ell I$ are invertible matrices for all integers $\ell\geq 0.$. If $r\leq s+2,$ for $\left\vert z\right\vert <1$%
%-----------------------------------------------------------------------------------------------------------------------------------------------------------------------------------
\begin{equation}
\begin{split}
&\;_{r}e_{s}\left(A_{1},A_{2},\ldots ,A_{r};B_{1},B_{2},\ldots ,B_{s};P,Q,z\right)\\
&=\Gamma ^{-1}\left( A_{i}\right) \Gamma ^{-1}\left( B_{j}-A_{i}\right)
\Gamma \left( B_{j}\right) \int_{0}^{1}t^{A_{i}-I}\left( 1-t\right)
^{B_{j}-A_{i}-I} \\
&\times\;_{r-1}e_{s-1}\left(
\begin{array}{c}
A_{1},\ldots ,A_{i-1},A_{i+1}\ldots ,A_{r}; \\
B_{1},\ldots ,B_{j-1},B_{j+1}\ldots ,B_{s}
\end{array}
;P,Q,zt\right)dt\label{3.38}
\end{split}
\end{equation}
%-----------------------------------------------------------------------------------------------------------------------------------------------------------------------------------
and
%-----------------------------------------------------------------------------------------------------------------------------------------------------------------------------------
\begin{equation}
\begin{split}
&\;_{r}E_{s}\left(A_{1},A_{2},\ldots ,A_{r};B_{1},B_{2},\ldots ,B_{s};P,Q,z\right)\\
&=\Gamma ^{-1}\left( A_{i}\right) \Gamma ^{-1}\left( B_{j}-A_{i}\right)
\Gamma \left( B_{j}\right) \int_{0}^{1}t^{A_{i}-I}\left( 1-t\right)
^{B_{j}-A_{i}-I} \\
&\times\;_{r-1}E_{s-1}\left(
\begin{array}{c}
A_{1},\ldots ,A_{i-1},A_{i+1}\ldots ,A_{r}; \\
B_{1},\ldots ,B_{j-1},B_{j+1}\ldots ,B_{s}
\end{array}
;P,Q,zt\right)dt.\label{3.39}
\end{split}
\end{equation}
%-----------------------------------------------------------------------------------------------------------------------------------------------------------------------------------
\end{theorem}
%-----------------------------------------------------------------------------------------------------------------------------------------------------------------------------------
\begin{proof}
%-----------------------------------------------------------------------------------------------------------------------------------------------------------------------------------
By generally for Eq. (\ref{3.26}), we have
%-----------------------------------------------------------------------------------------------------------------------------------------------------------------------------------
\begin{equation*}
\begin{split}
\left( A_{i}\right) _{\ell}[\left( B_{j}\right) _{\ell}]^{-1}=\Gamma ^{-1}\left(
A_{i}\right) \Gamma ^{-1}\left( B_{j}-A_{i}\right) \Gamma \left(
B_{j}\right) \int_{0}^{1}t^{A_{i}+\left( \ell-1\right) I}\left( 1-t\right)
^{B_{j}-A_{i}-I}dt
\end{split}
\end{equation*}
%-----------------------------------------------------------------------------------------------------------------------------------------------------------------------------------
where $A_{i}B_{j}=B_{j}A_{i}$. Also we have
%-----------------------------------------------------------------------------------------------------------------------------------------------------------------------------------
\begin{equation*}
\begin{split}
&\;_{r}e_{s}\left(
\begin{array}{c}
A_{1},A_{2},\ldots ,A_{r}; \\
B_{1},B_{2},\ldots ,B_{s};
\end{array}
z\right)  \\
&=\sum_{\ell=0}^{\infty }\frac{z^{\ell}}{k!}\left( A_{1}\right) _{\ell}\ldots \left(
A_{i}\right) _{\ell}\ldots \left( A_{r}\right) _{\ell}[\left( B_{1}\right)
_{\ell}]^{-1}\ldots [\left( B_{j}\right) _{\ell}]^{-1}\ldots [\left( B_{s}\right)
_{\ell}]^{-1}\\
&=\sum_{\ell=0}^{\infty }\frac{z^{\ell}}{k!}\left( A_{1}\right) _{\ell}\ldots \left(
A_{i-1}\right) _{\ell}\left( A_{i+1}\right) _{\ell}\ldots \left( A_{r}\right)
_{\ell}[\left( B_{1}\right) _{\ell}]^{-1}\ldots [\left( B_{j-1}\right)
_{\ell}]^{-1}[\left( B_{j+1}\right) _{\ell}]^{-1}
\\
&\ldots [\left( B_{s}\right)_{\ell}]^{-1}\times \Gamma ^{-1}\left( A_{i}\right) \Gamma ^{-1}\left(
B_{j}-A_{i}\right) \Gamma \left( B_{j}\right)\int_{0}^{1}t^{A_{i}+\left( n-1\right) I}\left(
1-t\right) ^{B_{j}-A_{i}-I}dt  \\
&=\Gamma ^{-1}\left( A_{i}\right) \Gamma ^{-1}\left( B_{j}-A_{i}\right)
\Gamma \left( B_{j}\right) \int_{0}^{1}t^{A_{i}-I}\left( 1-t\right)
^{B_{j}-A_{i}-I} \\
&\times \sum_{\ell=0}^{\infty }\frac{\left( zt\right) ^{\ell}}{k!}\left(
A_{1}\right) _{\ell}\ldots \left( A_{i-1}\right) _{\ell}\left( A_{i+1}\right)
_{\ell}\ldots \left( A_{r}\right) _{\ell}\\
&[\left( B_{1}\right) _{\ell}]^{-1}\ldots[\left( B_{j-1}\right) _{\ell}]^{-1}[\left( B_{j+1}\right) _{\ell}]^{-1}\ldots[\left(
B_{s}\right) _{\ell}]^{-1}dt \\
&=\Gamma ^{-1}\left( A_{i}\right) \Gamma ^{-1}\left( B_{j}-A_{i}\right)
\Gamma \left( B_{j}\right) \int_{0}^{1}t^{A_{i}-I}\left( 1-t\right)
^{B_{j}-A_{i}-I}\text{} \\
&\times\;_{r-1}e_{s-1}\left(
\begin{array}{c}
A_{1},\ldots, A_{i-1},A_{i+1},\ldots, A_{r}; \\
B_{1},\ldots, B_{j-1},B_{j+1},\ldots, B_{s};
\end{array}
zt\right)dt.
\end{split}
\end{equation*}
%-----------------------------------------------------------------------------------------------------------------------------------------------------------------------------------
%-----------------------------------------------------------------------------------------------------------------------------------------------------------------------------------
\end{proof}
%-----------------------------------------------------------------------------------------------------------------------------------------------------------------------------------
%---------------------------------------------------------------------------------------------------------------------------------------------------------------------------------------------
\begin{theorem}
%---------------------------------------------------------------------------------------------------------------------------------------------------------------------------------------------
The $\;_{r}e_{s}$ and $\;_{r}E_{s}$ matrix functions have the following integral representation
%---------------------------------------------------------------------------------------------------------------------------------------------------------------------------------------------
\begin{equation}
\begin{split}
&\;_{r}e_{s}(A_{1},A_{2},\ldots,A_{r};B_{1},B_{2},\ldots,B_{s};P,Q,z)=\Gamma^{-1}(A_{1})\\
&\int_{0}^{\infty}t^{A_{1}-I}e^{-t}\;_{r-1}e_{s}(A_{2},\ldots,A_{r};B_{1},B_{2},\ldots,B_{s};P,Q,zt)dt\label{3.40}
\end{split}
\end{equation}
%---------------------------------------------------------------------------------------------------------------------------------------------------------------------------------------------
and
%---------------------------------------------------------------------------------------------------------------------------------------------------------------------------------------------
\begin{equation}
\begin{split}
&\;_{r}E_{s}(A_{1},A_{2},\ldots,A_{r};B_{1},B_{2},\ldots,B_{s};P,Q,z)=\Gamma^{-1}(A_{1})\\
&\int_{0}^{\infty}t^{A_{1}-I}e^{-t}\;_{r-1}E_{s}(A_{2},\ldots,A_{r};B_{1},B_{2},\ldots,B_{s};P,Q,zt)dt.\label{3.41}
\end{split}
\end{equation}
%---------------------------------------------------------------------------------------------------------------------------------------------------------------------------------------------
%---------------------------------------------------------------------------------------------------------------------------------------------------------------------------------------------
\end{theorem}
%---------------------------------------------------------------------------------------------------------------------------------------------------------------------------------------------
%---------------------------------------------------------------------------------------------------------------------------------------------------------------------------------------------
\begin{proof}
%---------------------------------------------------------------------------------------------------------------------------------------------------------------------------------------------
Using the definition Gamma matrix function
%---------------------------------------------------------------------------------------------------------------------------------------------------------------------------------------------
\begin{equation*}
\begin{split}
\Gamma(A_{1}+\ell I)=\int_{0}^{\infty}e^{-t}t^{A_{1}+(\ell-1)I}dt,
\end{split}
\end{equation*}
%---------------------------------------------------------------------------------------------------------------------------------------------------------------------------------------------
we get (\ref{3.40}).
%---------------------------------------------------------------------------------------------------------------------------------------------------------------------------------------------
\end{proof}
%---------------------------------------------------------------------------------------------------------------------------------------------------------------------------------------------
\begin{theorem}
%---------------------------------------------------------------------------------------------------------------------------------------------------------------------------------------
We have the following recurrence matrix relation for the generalized incomplete exponential $\;_{r}E_{s}(x,P,Q;z)$ matrix functions:
%---------------------------------------------------------------------------------------------------------------------------------------------------------------------------------------
\begin{equation}
\begin{split}
(A_1-B_1+I) \ {\displaystyle _re_s\!\left[(x,P,Q;z)\left |{\begin{matrix} A_1,A_2,\ldots,A_{r}\\ B_1,B_{2},\ldots,B_{s}\end{matrix}} \right.\right]} &= A_{1} \ {\displaystyle _re_s\!\left[(x,P,Q;z)\left |{\begin{matrix} A_1+I,A_2,\ldots,A_{r}\\ B_1,B_{2},\ldots,B_{s}\end{matrix}} \right.\right]} \\ &- (B_1 -I)\ {\displaystyle _re_s\!\left[(x,P,Q;z)\left |{\begin{matrix} A_1,A_2,\ldots,A_{r}\\ B_1-I,B_{2},\ldots,B_{s}\end{matrix}} \right.\right]}.\label{3.42}
\end{split}
\end{equation}
%---------------------------------------------------------------------------------------------------------------------------------------------------------------------------------------
and
%---------------------------------------------------------------------------------------------------------------------------------------------------------------------------------------
\begin{equation}
\begin{split}
(A_1-B_1+I) \ {\displaystyle _rE_s\!\left[(x,P,Q;z)\left |{\begin{matrix} A_1,A_2,\ldots,A_{r}\\ B_1,B_{2},\ldots,B_{s}\end{matrix}} \right.\right]} &= A_{1} \ {\displaystyle _rE_s\!\left[(x,P,Q;z)\left |{\begin{matrix} A_1+I,A_2,\ldots,A_{r}\\ B_1,B_{2},\ldots,B_{s}\end{matrix}} \right.\right]} \\ &- (B_1 -I)\ {\displaystyle _rE_s\!\left[(x,P,Q;z)\left |{\begin{matrix} A_1,A_2,\ldots,A_{r}\\ B_1-I,B_{2},\ldots,B_{s}\end{matrix}} \right.\right]}.\label{3.43}
\end{split}
\end{equation}
%---------------------------------------------------------------------------------------------------------------------------------------------------------------------------------------------
%---------------------------------------------------------------------------------------------------------------------------------------------------------------------------------------------
\end{theorem}
%---------------------------------------------------------------------------------------------------------------------------------------------------------------------------------------
\begin{proof}
%---------------------------------------------------------------------------------------------------------------------------------------------------------------------------------------
Let $\mathfrak{R}$ be the l.h.s. of (\ref{3.42}). Then, using (\ref{3.7}), we get
%---------------------------------------------------------------------------------------------------------------------------------------------------------------------------------------
\begin{equation}
\begin{split}
\mathfrak{R}=  \sum_{\ell=0}^{\infty} \gamma(P \ell+Q,x)\Gamma^{-1}(P \ell+Q)A_1(A_1+I)_\ell(A_2)_\ell \ldots(A_r)_\ell[(B_1)_\ell]^{-1}\ldots [(B_s)_\ell]^{-1}\frac{z^\ell}{\ell !}\\
 - \sum_{\ell=0}^{\infty} \gamma(P \ell+Q,x)[\gamma(P \ell+Q)]^{-1}(A_1)_\ell(A_2)_\ell\ldots(A_r)_\ell(B_1-I)[(B_1-I)_\ell]^{-1}\ldots [(B_s)_\ell]^{-1}\frac{z^\ell}{\ell !}.\label{3.44}
\end{split}
\end{equation}
%---------------------------------------------------------------------------------------------------------------------------------------------------------------------------------------
%---------------------------------------------------------------------------------------------------------------------------------------------------------------------------------------
Using the relations
%---------------------------------------------------------------------------------------------------------------------------------------------------------------------------------------
%---------------------------------------------------------------------------------------------------------------------------------------------------------------------------------------
\begin{equation*}
\begin{split}
A_1 (A_1+I)_\ell &= (A_1+\ell I) (A_1)_\ell, \\
(B_1-I)[(B_1-I)_\ell]^{-1}&=(B_1+(\ell-1)I)[(B_1)_\ell]^{-1},
\end{split}
\end{equation*}
%---------------------------------------------------------------------------------------------------------------------------------------------------------------------------------------
this yields the left hand side of (\ref{3.42}).
%---------------------------------------------------------------------------------------------------------------------------------------------------------------------------------------
\end{proof}
%---------------------------------------------------------------------------------------------------------------------------------------------------------------------------------------
%---------------------------------------------------------------------------------------------------------------------------------------------------------------------------------------
We need the following result for proving Theorem \ref{thm3.17}.
%---------------------------------------------------------------------------------------------------------------------------------------------------------------------------------------
\begin{equation}
\begin{split}
(C+(1-k)I)_{n} =\Gamma(C+(1-k+n)I)\Gamma^{-1}(C+(1-k)I)\\
= (C+I)_{n} \cdot ((k-1)I-C)_{k}\cdot[((k-n-1)I-C)_{k}]^{-1}  , \;\; (k,n \in \mathbf{N}_0).\label{3.45}
\end{split}
\end{equation}
%---------------------------------------------------------------------------------------------------------------------------------------------------------------------------------------

%---------------------------------------------------------------------------------------------------------------------------------------------------------------------------------------
\begin{theorem}\label{thm3.17}
%---------------------------------------------------------------------------------------------------------------------------------------------------------------------------------------
For the generalized incomplete exponential $\;_{r}E_{s}(x,P,Q;z)$ and $\;_{r}e_{s}(x,P,Q;z)$ matrix functions, the following infinite summation relations holds true:
%---------------------------------------------------------------------------------------------------------------------------------------------------------------------------------------
%---------------------------------------------------------------------------------------------------------------------------------------------------------------------------------------
\begin{equation}
\begin{split}
& \sum_{k=0}^{\infty}  \frac{((k-1)I-C)_{k}}{k!} \ {\displaystyle _re_{s+1}\!\left[(x,P,Q;z)\left |{\begin{matrix} A_1,A_2,\dots,A_r\\ C-k+1, B_1,B_2,\dots,B_s\end{matrix}} \right.\right]} \ u^k = (1-u)^{C} \\ & \;\;\;\;\;\;\;\;\;\;\;\;\;\;\;\;\;\;\;\;\;\;\;\;\;\;\;\;\;\;\;\;\;\;\;\;\;\;\;\;\;\;\;\;\;\;\;\;\;\; \times {\displaystyle _re_{s+1}\!\left[(x,P,Q;z(1-u))\left |{\begin{matrix} A_1,A_2,\dots,A_r\\ I-C, B_1,B_2,\dots,B_s\end{matrix}} \right.\right]},\label{3.46}
\end{split}
\end{equation}
%---------------------------------------------------------------------------------------------------------------------------------------------------------------------------------------

%---------------------------------------------------------------------------------------------------------------------------------------------------------------------------------------
\begin{equation}
\begin{split}
& \sum_{k=0}^{\infty}  \frac{((k-1)I-C)_{k}}{k!} \ {\displaystyle _rE_{s+1}\!\left[(x,P,Q;z)\left |{\begin{matrix} A_1,A_2,\dots,A_r\\ C-k+1, B_1,B_2,\dots,B_s\end{matrix}} \right.\right]} \ u^k = (1-u)^{C} \\ & \;\;\;\;\;\;\;\;\;\;\;\;\;\;\;\;\;\;\;\;\;\;\;\;\;\;\;\;\;\;\;\;\;\;\;\;\;\;\;\;\;\;\;\;\;\;\;\; \times {\displaystyle _rE_{s+1}\!\left[(x,P,Q;v(1-u))\left |{\begin{matrix} A_1,A_2,\dots,A_r\\ I-C, B_1,B_2,\dots,B_s\end{matrix}} \right.\right]}, \label{3.47}
\end{split}
\end{equation}
%---------------------------------------------------------------------------------------------------------------------------------------------------------------------------------------
where $P,Q,C\in\mathbf{C}$ and $\left| u \right|<1$.
%---------------------------------------------------------------------------------------------------------------------------------------------------------------------------------------
\end{theorem}
%---------------------------------------------------------------------------------------------------------------------------------------------------------------------------------------
\begin{proof}
%---------------------------------------------------------------------------------------------------------------------------------------------------------------------------------------
Let $\mathcal{L}$ be the left hand side of (\ref{3.46}) and applying (\ref{3.6}), we arrive at
%---------------------------------------------------------------------------------------------------------------------------------------------------------------------------------------
\begin{equation*}
\begin{split}
\mathcal{L} =& \sum_{k=0}^{\infty}  \frac{((k-1)-C)_{k}}{k!} \sum_{n=0}^{\infty}  \gamma(P n+Q,x)\Gamma^{-1}(P n+Q)(A_1)_n(A_2)_n\dots(A_r)_n\\
&\times [(C+(1-k)I)_n]^{-1} [(B_1)_n]^{-1}[(B_2)_n]^{-1}\dots[(B_s)_n]^{-1} \frac{z^n}{n!} t^k.
\end{split}
\end{equation*}
%---------------------------------------------------------------------------------------------------------------------------------------------------------------------------------------
By reversing the order of summation and using the result defined in (\ref{3.51}), we obtain
%---------------------------------------------------------------------------------------------------------------------------------------------------------------------------------------
\begin{equation}
\begin{split}
\mathcal{C}_1 =& \sum_{n=0}^{\infty} \gamma(P n+Q,x)\Gamma^{-1}(P n+Q)(A_1)_n(A_2)_n\dots(A_r)_n[(C+I)_n]^{-1} \\
&\times [(B_1)_n]^{-1}[(B_2)_n]^{-1}\dots[(B_s)_n]^{-1} \frac{z^n}{n!} \cdot \sum_{k=0}^{\infty} \frac{((k-n-1)I-C)_{k}}{k!} t^k.\label{3.48}
\end{split}
\end{equation}
%---------------------------------------------------------------------------------------------------------------------------------------------------------------------------------------
Now we find the inner sum in (\ref{3.48}) by using binomial expansion
%---------------------------------------------------------------------------------------------------------------------------------------------------------------------------------------
\begin{equation}
\begin{split}
\sum_{k=0}^{\infty} \frac{((k-1)I-C)_{k}}{k!} t^k = \sum_{k=0}^{\infty} \frac{(-C)_{k}}{k!} t^k =(1-t)^C, \; \; \left| t \right|<1.\label{3.49}
\end{split}
\end{equation}
%---------------------------------------------------------------------------------------------------------------------------------------------------------------------------------------
Replacing the inner sum of (\ref{3.48}) and using (\ref{3.49}), this yields the r.h.s. of (\ref{3.46}). Similarly, one can prove (\ref{3.55}) as asserted by Theorem \ref{thm3.17}.
%---------------------------------------------------------------------------------------------------------------------------------------------------------------------------------------
\end{proof}
%---------------------------------------------------------------------------------------------------------------------------------------------------------------------------------------
\begin{remark}
%---------------------------------------------------------------------------------------------------------------------------------------------------------------------------------------
If we add infinite summation formulas (\ref{3.46}) and (\ref{3.47}) and using (\ref{3.49}), then we have the following interesting infinite summation formula
%---------------------------------------------------------------------------------------------------------------------------------------------------------------------------------------
%---------------------------------------------------------------------------------------------------------------------------------------------------------------------------------------
\begin{equation}
\begin{split}
& \sum_{k=0}^{\infty} \frac{(k-1)I-C)_{k}}{k!} \ {\displaystyle _rF_{s+1}\!\left[z\left |{\begin{matrix} A_1,A_2,\dots,A_r\\ C+(1-k)I, B_1,B_2,\dots,B_s\end{matrix}} \right.\right]} \ t^k = (1-t)^{C} \\ & \;\;\;\;\;\;\;\;\;\;\;\;\;\;\;\;\;\;\;\;\;\;\;\;\;\;\;\;\;\;\;\;\;\;\;\;\;\;\;\;\;\;\;\;\;\;\;\;\;\;\;\;\;\; \times {\displaystyle _rF_{s+1}\!\left[z(1-t)\left |{\begin{matrix} A_1,A_2,\dots,A_r\\ I-C, B_1,B_2,\dots,B_s\end{matrix}} \right.\right]},\label{3.50}
\end{split}
\end{equation}
%---------------------------------------------------------------------------------------------------------------------------------------------------------------------------------------
where $C\in\mathbf{C}$ and $\left| t \right|<1$.
%---------------------------------------------------------------------------------------------------------------------------------------------------------------------------------------
\end{remark}
%---------------------------------------------------------------------------------------------------------------------------------------------------------------------------------------
\begin{definition}
%---------------------------------------------------------------------------------------------------------------------------------------------------------------------------------------
For the extension of the generalized incomplete exponential matrix functions (\ref{3.6}) and (\ref{3.7}), we consider two sequences $\{\Psi _p^{\left( {C,\kappa } \right)}\left( z \right)\}_{p=0}^\infty$ and $\{\Phi _p^{\left( {C,\kappa } \right)}\left( z \right)\}_{p=0}^\infty$ as,
%---------------------------------------------------------------------------------------------------------------------------------------------------------------------------------------
\begin{equation}
\begin{split}
\Psi _p^{\left( {C,\kappa } \right)}\left( z \right) &=  {\displaystyle \Psi _p^{\left( {C,\kappa } \right)} \!\left[(x,P,Q;z)\left |{\begin{matrix} A_1,A_2,\dots,A_r\\ B_1,B_2,\dots,B_s\end{matrix}} \right.\right]},\\
&= {\displaystyle _re_{s+\kappa}\!\left[(x,P,Q;z)\left |{\begin{matrix} A_1,A_2,\dots,A_r\\\Delta(\kappa, (1-p)I-C), B_1,B_2,\dots,B_s\end{matrix}} \right.\right]} \label{3.51}
\end{split}
\end{equation}
%---------------------------------------------------------------------------------------------------------------------------------------------------------------------------------------
%---------------------------------------------------------------------------------------------------------------------------------------------------------------------------------------
and
%---------------------------------------------------------------------------------------------------------------------------------------------------------------------------------------
\begin{equation}
\begin{split}
\Phi _p^{\left( {C,\kappa } \right)}\left( z \right) &=  {\displaystyle \Phi _p^{\left( {C,\kappa } \right)} \!\left[(x,P,Q;z)\left |{\begin{matrix} A_1,A_2,\dots,A_r\\ B_1,B_2,\dots,B_s\end{matrix}} \right.\right]},\\
&= {\displaystyle _rE_{s+\kappa}\!\left[(x,P,Q;z)\left |{\begin{matrix} A_1,A_2,\dots,A_r\\\Delta(\kappa, (1-p)I-C), B_1,B_2,\dots,B_s\end{matrix}} \right.\right]},\label{3.52}
\end{split}
\end{equation}
%---------------------------------------------------------------------------------------------------------------------------------------------------------------------------------------
where $\Delta(\kappa, A)$ abbreviates the array of $\kappa$ parameters as follows:
%---------------------------------------------------------------------------------------------------------------------------------------------------------------------------------------
\begin{equation*}
\begin{split}
\frac{A}{\kappa}, \frac{A+I}{\kappa}, \frac{A+2I}{\kappa}, \dots ,\frac{A+(\kappa-1)I}{\kappa}, \;\; (A\in\mathbf{C}, \ \kappa\in\mathbf{N}).
\end{split}
\end{equation*}
%---------------------------------------------------------------------------------------------------------------------------------------------------------------------------------------
%---------------------------------------------------------------------------------------------------------------------------------------------------------------------------------------
\end{definition}
%---------------------------------------------------------------------------------------------------------------------------------------------------------------------------------------
Now, we derive the generalization of the infinite summation formulas for the generalized incomplete exponential matrix functions (\ref{3.6}) and (\ref{3.7}) asserted in the following theorem.
%---------------------------------------------------------------------------------------------------------------------------------------------------------------------------------------
\begin{theorem}\label{thm3.18}
%---------------------------------------------------------------------------------------------------------------------------------------------------------------------------------------
The following infinite summation formulas holds true for the following two sequences $\{\Psi _p^{\left( {C,\kappa } \right)}\left( z \right)\}_{p=0}^\infty$ and $\{\Phi _p^{\left( {C,\kappa } \right)}\left( z \right)\}_{p=0}^\infty$ :
%---------------------------------------------------------------------------------------------------------------------------------------------------------------------------------------
\begin{equation}
\begin{split}
\sum_{p=0}^{\infty}  \frac{(C+(m+p-1)I)_{r}}{p!}   \Psi _{m+p}^{\left( {C,\kappa } \right)}\left( z \right) t^p = (1-t)^{-C-mI} \ \Psi _{m}^{\left( {C,\kappa } \right)}\left( z(1-t)^\kappa \right),\label{3.53}
\end{split}
\end{equation}
%---------------------------------------------------------------------------------------------------------------------------------------------------------------------------------------
\begin{equation}
\begin{split}
\sum_{p=0}^{\infty}  \frac{(C+(m+p-1)I)_{r}}{p!} \  \Phi _{m+p}^{\left( {C,\kappa } \right)}\left( z \right) t^p = (1-t)^{-C-mI} \ \Phi _{m}^{\left( {C,\kappa } \right)}\left( z(1-t)^\kappa \right),\label{3.54}
\end{split}
\end{equation}
%---------------------------------------------------------------------------------------------------------------------------------------------------------------------------------------
where $x\geq0; m\in\mathbf{N}_0;C\in\mathbf{C}, \ \kappa\in\mathbf{N}; \ \left| u \right|<1.$
%---------------------------------------------------------------------------------------------------------------------------------------------------------------------------------------
\end{theorem}
%---------------------------------------------------------------------------------------------------------------------------------------------------------------------------------------
\begin{proof}
%---------------------------------------------------------------------------------------------------------------------------------------------------------------------------------------
Proof of Theorem \ref{thm3.18} is similar to Theorem \ref{thm3.17}. In aforesaid theorem, one can use the following identity,
%---------------------------------------------------------------------------------------------------------------------------------------------------------------------------------------
\begin{align}
(I-C-(m+p)I)_{\kappa \ell} = ((1-m)I-C)_{\kappa \ell}(C+(m+p-1)I)_{p} [(C+(m-\kappa \ell+p-1)I)_{p}]^{-1}\label{3.55}
\end{align}
%---------------------------------------------------------------------------------------------------------------------------------------------------------------------------------------
where $p,m,\ell\in\mathbf{N}_0;\kappa\in\mathbf{N}$.
%---------------------------------------------------------------------------------------------------------------------------------------------------------------------------------------
\end{proof}
%---------------------------------------------------------------------------------------------------------------------------------------------------------------------------------------
\begin{remark}
%---------------------------------------------------------------------------------------------------------------------------------------------------------------------------------------
On setting $\kappa=1$ and replacing $C$ by $C-mI$ in (\ref{3.53}), this can easily reduce to the result (\ref{3.46}).
%---------------------------------------------------------------------------------------------------------------------------------------------------------------------------------------
\end{remark}
%---------------------------------------------------------------------------------------------------------------------------------------------------------------------------------------------
%  \section{Applications}
%---------------------------------------------------------------------------------------------------------------------------------------------------------------------------------------------

%---------------------------------------------------------------------------------------------------------------------------------------------------------------------------------------------
\section{Conclusion}
%---------------------------------------------------------------------------------------------------------------------------------------------------------------------------------------------
This paper introduces a new method for investigating some interesting results for incomplete gamma matrix functions $\gamma(Q, x)$ and $\Gamma(Q, x)$, including  recurrence matrix relations, differentiation formulas,  as well as analytical and fractional integral properties of incomplete gamma matrix functions. Additionally, some relevant characteristics of these functions such as integral representations functions have been established. We derive the generalization of the infinite summation formulas for the generalized incomplete exponential matrix functions  along with the generalized incomplete exponential matrix functions with the integral representation. Additionally, we derive the addition formula for addition of two arguments, multiplication formula for multiplication of two arguments, and  some concluding remarks and consequences of the results have also been established. In mathematical physics, fractional calculus, and classical analysis, the results of this work can be extremely important.
%---------------------------------------------------------------------------------------------------------------------------------------------------------------------------------------------
%---------------------------------------------------------------------------------------------------------------------------------------------------------------------------------------------
%---------------------------------------------------------------------------------------------------------------------------------------------------------------------------------------------
%---------------------------------------------------------------------------------------------------------------------------------------------------------------------------------------------

%---------------------------------------------------------------------------------------------------------------------------------------------------------------------------------------------
\end{document}